\documentclass[12pt,thmsa]{article}
\usepackage{amsmath}
\usepackage{amssymb}
\usepackage{graphicx}
\usepackage{url} % not crucial - just used below for the URL 
\usepackage{stmaryrd}
\usepackage{mathrsfs}
\usepackage{algorithm}
\usepackage{algorithmic}
\usepackage{epstopdf}% To incorporate .eps illustrations using PDFLaTeX, etc.
\usepackage[caption=false]{subfig}% Support for small, `sub' figures and tables

\newtheorem{remark}{Remark}

\begin{document}

\begin{center}
{\Large \textbf{ A method for variable selection in a multivariate functional linear  regression model}}

\bigskip

Alban  MBINA MBINA and  Guy Martial  NKIET 

\bigskip

\textsuperscript{}URMI, Universit\'{e} des Sciences et Techniques de Masuku,  Franceville, Gabon.

\bigskip

E-mail adresses : alban.mbinambina@univ-masuku.org,    guymartial.nkiet@univ-masuku.org.

\bigskip
\end{center}

\noindent\textbf{Abstract.}  We propose a new
variable selection procedure for a functional
linear model with  multiple scalar responses  and multiple functional predictors. This method is based on basis expansions of the involved functional predictors and coefficients that  lead to a multivariate linear regression model.  Then a  criterion  by means of which   the variable selection problem reduces to that of estimating a suitable set is introduced. Estimation of this set is achieved  by using appropriate penalizations of estimates of this  criterion, so leading to our proposal. A simulation study that permits to investigate the effectiveness of the proposed approach and to compare it with  
existing methods is given. 
\bigskip

\noindent\textbf{AMS 1991 subject classifications: }62H99, 62J05.

\noindent\textbf{Key words: }Variable selection; Functional linear; Selection criterion; Functional Data Analysis.
\section{Introduction}\label{sec1}

In statistical modeling, an  usual  approach that consists   in  determining a model linking a response variable to a set of predictor variables has given rise to numerous regression models, including the  multivariate  linear regression model which has been intensively studied from different aspects for many years. One of the most crucial issues  related to this model is the variable selection problem which arises when one  has  to determine, among a number of predictors  which can be large, the variables which are really relevant to explain the response. So, many methods offering solutions for this problem have been proposed in the literature. Surveys on earlier works in this field can be found  in   Hocking (1976), Thomson (1978a), Thomson (1978b),  and some  recent references on this topic are  Ranciati et al. (2019), Bizuayehu et al. (2022),   Mbina Mbina et al. (2023)  and Wei and Yu (2023). On the other hand, statistical methods for processing data in the form of curves have had significant development over the last twenty years, thus allowing the emergence of a very active field of statistics called  Functional Data Analysis (FDA) which  has received considerable attention due to its  large number of applications (see, e.g.,  Ramsay and Silverman (2005), Ferraty and Vieu (2006), Horv\'ath and   Kokoszka (2012), Kokoszka and   Reimherr (2017)). Surveys on recent  developments on FDA can be found in Goia16 and Vieu (2016) and Aneiros et al. (2019). In the context of FDA also, linear regression models have been introduced with the aim of describing the relationships between several functional variables and one or more response variables which may also   be of functional nature or not. The functional linear model, where there is only one predictor of functional nature, was first considered (see, e.g., Cardot et al. (1999), Cardot et al. (2003)). Later, generalizations of this model including several functional predictors, prompted by applications, were then addressed. Here too, the problem of selecting the  functional predictors that are really  relevant to explain the response variable is of  great  importance for modeling. However, only a few authors have considered variable selection in functional regression analysis with several functional predictors; a survey on works in  this field can be found in Aneiros et al. (2022). Matsui and Konishi (2011)  adapted an approach introduced by Fan and Li (2001) to obtain a method based on  $L_1$ regularization. Lian(201)  studied selection of relevant functional variables by using functional principal components basis expansions. Collazos et al. (2016)  proposed a method based on testing for the nullity of functional coefficients:  a covariate is dropped  from the model  when   the null hypothesis that its corresponding
parameter  is equal to zero is not rejected. For doing that, they introduced a likelihood ratio type test, where restricted and full models are estimated
through the B-Splines basis expansions of both coefficients and functional predictors.  Liu et al. (2018)   proposed a
functional variable selection procedure  using the technique of Gram–Schmidt orthogonalization
to remove the irrelevant predictors. Smaga and Matsui (2018)  introduced two methods based on random subspace method of Mielniczuk and Teisseyre (2014). Matsui and Umezu (2020) considered  the use of
sparse regularization in the construction of a functional regression model with
functional predictors and multiple scalar responses. The aforementioned  methods  are essentially based on extensions to the functional case of methods suitable for multiple linear regression. This extension is made possible thanks to basis expansions   of the involved functional  variables and   coefficients as tackled, for example, in   G\'orecki e al. (2018).

In this paper, we propose a new method for variable selection in multivariate functional linear regression model by extending the approach of Mbina Mbina et al. (2023) from basis expansions allowing to transform the initial model to a multivariate linear regression model as   described in Section \ref{sec2}. This method is based on   a  criterion, introduced in   Section \ref{sec2}, by means of which   the variable selection problem reduces to that of estimating a suitable set. Then our proposal for variable selection is achieved in  Section  \ref{sec:MLE} from an estimate of  this set obtained by using appropriate penalizations of estimates of the aforementioned criterion.
The effectiveness of the proposed variable selection strategy is investigated in Section  \ref{sec4:sim} through Monte Carlo simulations which set up   comparison with a random subspace method of Smaga and Matsui (2018)  and the group SCAD method of Matsui and Konishi (2011).

\section{MFLR  model and variable selection problem}
\label{sec2}
\noindent

\noindent In this  section we first define the multivariate functional linear regression (MFLR)  model that is  used, then a transformation of this model  from  basis representations  of the involved functional variables and coeffcients is obtained. Finally, the criterion that is used for dealing with variable selection is specified.

\subsection{MFLR model}
For $(p,q)\in\left(\mathbb{N}^\ast \right)^2$ and $\ell=1,\cdots,p$, $j=1,\cdots,q$, we consider   real random variables $Y_j$  and processes $\{X_\ell(t);\,t\in\mathcal{I}_\ell\}$, where  $\mathcal{I}_\ell$ is an interval of $\mathbb{R}$. Assuming that  $X_\ell\in L^2(\mathcal{I}_\ell)$, we deal with the MFLR  model given by:
\begin{equation}\label{model}
Y_j = \sum_{\ell=1}^p\int_{\mathcal{I}_\ell}B_{j\ell}(t)\,X_{\ell}(t) \,\,dt + \varepsilon_j
\end{equation}
where the $B_{j\ell}$s  are functional coefficients  and $\varepsilon=\left(\varepsilon_1,\cdots,\varepsilon_q\right)^T$  is a random vector with values in  $\mathbb{R}^{q}$ with mean $0$ and unknown covariance matrix,  and which is independent of   $X=(X_1,\cdots,X_p)$. We are interested in variable selection in model (\ref{model}), that is determining the explanatory functional variables whose functional coefficients  are not null functions,   from an i.i.d. sample $\left\{\left(Y^{(i)},X^{(i)}\right)\right\}_{1\leq i\leq n}$  of  $\left(Y,X\right)$, where: 
\begin{equation}\label{sample}
Y=\left(Y_1,\cdots,Y_q\right)^T,\,\,\,Y^{(i)}=\left(Y_1^{(i)},\cdots,Y_q^{(i)}\right)^T\,\,\,\textrm{ and }\,\,\,X^{(i)}=\left(X_1^{(i)},\cdots,X_p^{(i)}\right).
\end{equation}
 Putting $\llbracket 1,m\rrbracket=\{1,\cdots,m\}$, we assume that the set
\begin{equation}\label{i0}
I_0=\{\ell\in \llbracket 1,p\rrbracket\,/\,B_{j\ell}(t)=0,\,\forall j\in\llbracket 1,q\rrbracket,\,\forall t\in\mathcal{I}_\ell\}
\end{equation}
is not empty, and we tackle the variable selection problem as a problem of estimating the set
$I_1=\llbracket 1,p\rrbracket-I_0$  containing the   integers $\ell\in \llbracket 1,p\rrbracket$ such that there exist $j\in\llbracket 1,q\rrbracket$ and            $t\in\mathcal{I}_\ell$ satisfying $B_{j\ell}(t)\neq 0$.

\subsection{Basis representation} 
Following \cite{gorecki18} we will simplify model (\ref{model}) by using basis representations of the functions involved in this model, so obtaining a multivariate linear regression model which will be  considered for variable selection purpose. For $\ell\in  \llbracket 1,p\rrbracket$, letting   $\left\{\phi_{k \ell}\right\}_{k\geq 1} $ be a basis of    $L^{2}(\mathcal{I}_\ell)$, we consider the following truncated representations
\begin{equation}\label{truncated}
B_{j\ell}(t)\simeq \sum_{k=1}^{d_\ell} \textrm{\textbf{b}}_{jk \ell}\phi_{k \ell}(t) \;\;\;  \textrm{and} \;\;\; X_{\ell}(t) \simeq \sum_{k=1}^{d_\ell} \textrm{\textbf{X}}_{k \ell}\phi_{k \ell}(t), %\;\; \Forall \;\; k=1, \cdots ,p, \;\; i=1, \cdots ,n.
\end{equation}
where  $d_\ell$ is a dimensionality parameter which is to be chosen from the above introduced sample $\left\{ X^{(i)}\right\}_{1\leq i\leq n}$  by using the Bayesian information criterion (BIC) as suggested in \cite{gorecki18}. More precisely,  if the $X_\ell^{(i)}$s are observed on a fine grid of points $t_1^{(\ell)},\cdots, t_{N_\ell}^{(\ell)}$ of $\mathcal{I}_\ell$, we chose the dimensions $d^{(i)}_\ell$ that minimize the BIC criterion given by
\begin{equation*}\label{bic}
\textrm{BIC}(i,\ell)=\ln\bigg(\sum_{r=1}^N\bigg(X^{(i)}_\ell(t^{(\ell)}_r)- \sum_{k=1}^{d^{(i)}_\ell} \textrm{\textbf{X}}^{(i)}_{k \ell}\phi_{k \ell}(t_r^{(\ell)})\bigg)^2\bigg)+(d^{(i)}_\ell+1)\bigg(\frac{\ln (N_\ell)}{N_\ell}\bigg),
\end{equation*}
where $\textrm{\textbf{X}}^{(i)}_{k \ell}$ is the coordinate  of $X_\ell^{(i)}$ related to $\phi_{k\ell}$, 
 then we take $d_\ell=\max\{d^{(1)}_\ell,\cdots,d^{(n)}_\ell\}$. The basis functions $\phi_{k \ell}$  may be, for instance, the spline or  Gaussian radial bases,   Fourier basis or   wavelet basis, depending on the nature of the data (e.g., Smaga and Matsui (2018)).  From equation \eqref{truncated} we have:
\begin{equation}\label{vectorapprox}
B_{j\ell}(t) = \textrm{\textbf{b}}_{j\centerdot\ell}^T \phi_{\centerdot\ell}(t)\;\;\; \textrm{and} \;\;\; X_{\ell}(t) = \textrm{\textbf{X}}_{\centerdot\ell}^T\phi_{\centerdot\ell}(t),
\end{equation}
where $ \textrm{\textbf{b}}_{j\centerdot\ell}=(\textrm{\textbf{b}}_{j1\ell}, \cdots ,\textrm{\textbf{b}}_{jd_\ell\ell})^T$,  $\textrm{\textbf{X}}_{\centerdot\ell}=(\textrm{\textbf{X}}_{1\ell}, \cdots ,\textrm{\textbf{X}}_{d_\ell\ell})^T$   and $$\phi_{\centerdot\ell}(t)=(\phi_{1\ell}(t), \cdots ,\phi_{d_\ell\ell}(t))^T.$$  Then, model \eqref{model} becomes:
\begin{equation}\label{model2}
Y_j = \sum_{\ell=1}^p\textrm{\textbf{b}}_{j\centerdot\ell}^T\textrm{\textbf{G}}_{\ell}\, \textrm{\textbf{X}}_{\centerdot\ell}+ \varepsilon_j,
\end{equation}
where $\textrm{\textbf{G}}_{\ell}$ is the $d_\ell\times d_\ell$ Gram matrix $\textrm{\textbf{G}}_{\ell}=\big(g_{\ell;km}\big)_{1\leq km\leq d}$, where
\begin{equation}\label{gramcoeff}
g_{\ell;km}=\int_{\mathcal{I}_\ell}\phi_{k\ell}(t)\,\phi_{m\ell}(t)\,\,dt.
\end{equation}
Clearly, \eqref{model2} can be writen as
\begin{equation}\label{modelfinal}
Y=\sum_{\ell=1}^pb_{\centerdot\centerdot\ell}^T\textrm{\textbf{G}}_{\ell}\, \textrm{\textbf{X}}_{\centerdot\ell}+ \varepsilon ,
\end{equation}
where $b_{\centerdot\centerdot\ell}$ is the $d_\ell\times q$ matrix obtained by stacking the $ \textrm{\textbf{b}}_{j\centerdot\ell}$s  into columns  as
\[
b_{\centerdot\centerdot\ell}=\left(
\begin{array}{ccccccc}
\textrm{\textbf{b}}_{1\centerdot\ell} &\vert&\textrm{\textbf{b}}_{2\centerdot\ell} &\vert&\cdots &\vert &\textrm{\textbf{b}}_{q\centerdot\ell}
\end{array}
\right).
\]
Then, from \eqref{i0} and  \eqref{vectorapprox} it is seen that the aforementioned variable selection problem leads to consider the   subset  of $  \llbracket 1,p\rrbracket$ given by
\[
\widetilde{I}_0=\{\ell\in  \llbracket 1,p\rrbracket/\,b_{\centerdot\centerdot\ell}=0\}
\]
and  to estimate  the set
\begin{equation*}\label{i1}
\widetilde{I}_1= \llbracket 1,p\rrbracket-\widetilde{I}_0=\{\ell\in \llbracket 1,p\rrbracket\,\, /\,\,b_{\centerdot\centerdot\ell}\neq0\}.
\end{equation*}
\begin{remark} 
The Gram matrix $\textrm{\textbf{G}}_{\ell}$ has different expressions depending on the chosen basis:

\noindent 1) If  $\left\{\phi_{k \ell}\right\}_{k\geq 1} $ is an orthonormal basis (e.g., Fourier basis or wavelets), then $\textrm{\textbf{G}}_{\ell}$ equals the $d_\ell\times d_\ell$ identity matrix, that is
\begin{equation}\label{gfourier}
g_{\ell;km}=\left\{
\begin{array}{lcl}
1 & &\textrm{if }k=m\\
 & & \\
0 & &\textrm{if }k\neq m
\end{array}
\right. .
\end{equation}

\noindent 2) If  $\left\{\phi_{k \ell}\right\}_{k\geq 1} $ is  a spline basis, then $\textrm{\textbf{G}}_{\ell}$ can be computed from \eqref{gramcoeff} by using the trapezoidal rule, that is
\begin{equation}\label{gspline}
g_{\ell;km}\simeq\frac{1}{2}\sum_{r=1}^{N-1}\left(t_{ r+1}^{(\ell)}-t_r^{(\ell)}\right)\, \left(\phi_{k\ell}(t_r^{(\ell)})\phi_{m\ell}(t_{r}^{(\ell)})+ \phi_{k\ell}(t_{r+1}^{(\ell)})\phi_{m\ell}(t_{r+1}^{(\ell)})\right).
\end{equation}

\noindent 3) If  Gaussian basis functions are considered, that is
\[
\phi_{k \ell}(t)=\exp\left\{-\frac{(t-c_{k\ell})^2}{2\gamma_{\ell}\,\sigma^2_{k\ell}}\right\},
\]
then (see Matsui et al. (2008))
\begin{equation}\label{ggaussian}
g_{\ell;km}=\frac{\sqrt{2\pi}\sigma_{k\ell}\sigma_{m\ell}}{\sqrt{\sigma_{k\ell}^2+\sigma_{m\ell}^2}}\exp\left\{-\frac{(c_{k\ell}-c_{m\ell})^2}{2\gamma_{\ell}\left(\sigma^2_{k\ell}+\sigma^2_{m\ell}\right)}\right\}.
\end{equation}

\end{remark}
\subsection{Criterion for variable selection}
In order to simplify the estimation of $\widetilde{I}_1$ we will first characterize this subset  by means of a criterion which introduced below. Considering the random vector $\mathcal{X}$ with values in $\mathbb{R}^{d}$, where $d=\sum_{\ell=1}^pd_\ell$,  defined as
\[
\mathcal{X}= \left(
 \begin{array}{c}
\textrm{\textbf{G}}_{1}\, \textrm{\textbf{X}}_{\centerdot 1}\\
\hline 
 \textrm{\textbf{G}}_{2}\, \textrm{\textbf{X}}_{\centerdot 2} \\
\hline 
 \vdots \\
\hline 
\textrm{\textbf{G}}_{p}\, \textrm{\textbf{X}}_{\centerdot p}
 \end{array}
\right) ,
\]
we  assume  that  $\mathbb{E}(\|Y\|_{\mathbb{R}^{q}}^{2}) < +\infty$ and $\mathbb{E}(\|\textrm{\textbf{X}}_{\centerdot \ell}\|_{\mathbb{R}^{d_\ell}}^{2}) < +\infty$ for any $\ell\in \llbracket 1,p\rrbracket$, where   $\|.\|_{\mathbb{R}^{m}}$ denotes the usual Euclidean norm of $\mathbb{R}^{m}$.   Then,  putting 
\[
\mu=\mathbb{E}\left( \mathcal{X}\right)= \left(
 \begin{array}{c}
\textrm{\textbf{G}}_{1}\, \textrm{\textbf{m}}_{\centerdot 1}\\
\hline 
 \textrm{\textbf{G}}_{2}\, \textrm{\textbf{m}}_{\centerdot 2} \\
\hline 
 \vdots \\
\hline 
\textrm{\textbf{G}}_{p}\, \textrm{\textbf{m}}_{\centerdot p}
 \end{array}
\right)\,\,\,\textrm{ and }\eta=\mathbb{E}\left( Y\right)= \left(
 \begin{array}{c}
\eta_1\\
 \vdots \\
\eta_q
 \end{array}
\right),
\]
where   $\textrm{\textbf{m}}_{\centerdot i}= \mathbb{E}\left(\textrm{\textbf{X}}_{\centerdot i}\right)$ and $\eta_j= \mathbb{E}\left(Y_j\right)$,   it is possible to consider   the $d\times d$ and $d\times q$  covariance and cross-covariance  matrices    given  by

\begin{equation*}\label{covop}
\mathcal{C}_{1} = \mathbb{E}\bigg(\left(\mathcal{X}-\mu\right)\left(\mathcal{X}-\mu\right)^{T}\bigg) = \left(
 \begin{array}{cccc}
  \textrm{\textbf{G}}_{1}V_{11} \textrm{\textbf{G}}_{1}^T&\textrm{\textbf{G}}_{1}V_{12} \textrm{\textbf{G}}_{2}^T &\cdots &\textrm{\textbf{G}}_{1}V_{1p} \textrm{\textbf{G}}_{p}^T\\
 \textrm{\textbf{G}}_{2}V_{21} \textrm{\textbf{G}}_{1}^T&\textrm{\textbf{G}}_{2}V_{22} \textrm{\textbf{G}}_{2}^T &\cdots &\textrm{\textbf{G}}_{2}V_{2p} \textrm{\textbf{G}}_{p}^T\\\
  \vdots &\vdots &\cdots &\vdots\\
\textrm{\textbf{G}}_{p}V_{p1} \textrm{\textbf{G}}_{1}^T&\textrm{\textbf{G}}_{p}V_{p2} \textrm{\textbf{G}}_{2}^T &\cdots &\textrm{\textbf{G}}_{p}V_{pp} \textrm{\textbf{G}}_{p}^T
 \end{array}
\right)
\end{equation*}
and
\begin{equation*}\label{covop2}
\mathcal{C}_{12}= \mathbb{E}\bigg(\left(\mathcal{X}-\mu\right)\left(Y-\eta\right)^{T}\bigg) = \left(
 \begin{array}{cccc}
  \textrm{\textbf{G}}_{1}W_{11}&\textrm{\textbf{G}}_{1}W_{12}  &\cdots &\textrm{\textbf{G}}_{1}W_{1q}  \\
 \textrm{\textbf{G}}_{2}W_{21} &\textrm{\textbf{G}}_{2}W_{22}  &\cdots &\textrm{\textbf{G}}_{2}W_{2q}\\\
  \vdots &\vdots &\cdots &\vdots\\
\textrm{\textbf{G}}_{p}W_{p1} &\textrm{\textbf{G}}_{p}W_{p2}  &\cdots &\textrm{\textbf{G}}_{p}W_{pq}
 \end{array}
\right),
\end{equation*}
where $V_{ij} = \mathbb{E}\left(\left( \textrm{\textbf{X}}_{\centerdot i}- \textrm{\textbf{m}}_{\centerdot i}\right)\left( \textrm{\textbf{X}}_{\centerdot j}- \textrm{\textbf{m}}_{\centerdot j}\right)^{T}\right)$ and  $W_{ij} = \mathbb{E}\left(\left( Y_j-\eta_j\right)\left(\textrm{\textbf{X}}_{\centerdot i} - \textrm{\textbf{m}}_{\centerdot i}\right)\right)$. Furthermore,  for a given subset $K=\{i_1,\cdots,i_k\}$ of $\llbracket 1,p\rrbracket$,  let us consider the  $(\sum_{r=1}^kd_{i_r})\times d$  matrix:
\[
A_K=\left(
\begin{array}{cccc}
a_{11}^{(K)} & a_{12}^{(K)} &\cdots &a_{1p}^{(K)}\\
a_{21}^{(K)} & a_{22}^{(K)} &\cdots &a_{2p}^{(K)}\\
\vdots & \vdots & \vdots &\vdots \\
a_{k1}^{(K)} & a_{k2}^{(K)} &\cdots &a_{kp}^{(K)}\\
\end{array}
\right)
\]
where, denoting by $\mathbb{I}_m$ the $m\times m$ identity matrix, we have
\[
a_{r j}^{(K)}=\left\{
\begin{array}{ccc}
\mathbb{I}_{d_j} &  &\textrm{ if }j=i_r\\
0  & &\textrm{ if }j\neq i_r\\
\end{array}
\right. , \,\,1\leq r\leq k,\,1\leq \ell\leq d_j.
\]
Then, putting $\Pi_{K}:=A_{K}^T\left(  A_{K}\mathcal{C}_{1}A_{K}^T\right)  ^{-1}A_{K}$, where $A^{-1}$ denotes the   inverse of the matrix $A$, we introduce the criterion
\begin{equation*}\label{crit}
\xi_{K} = \| \mathcal{C}_{12} -\mathcal{C}_{1}  \Pi_{K} \mathcal{C}_{12}  \|,
\end{equation*}
where $\Vert\cdot\Vert$ denotes the usual matrices   norm given by $\Vert A\Vert^2=\textrm{tr}\left(A\, A^T\right)$. This criterion it of a type which was considered in Mbina Mbina et al. (2023)  where it is  shown  that it measures  a  distance between the matrix of coefficients of model \eqref{modelfinal}  when the whole variables are considered and the matrix of  coefficients of the same model when only the variables whose indices belong to $K$ are considered, so measuring the relevance of these later variables for variable selection.  Using  this criterion we get  a more explicit expression of $\widetilde{I}_1$. Indeed,  similarly than in Mbina Mbina et al. (2023), $K$ is included in  $\widetilde{I}_1$ if and only if  $\xi_{K} = 0$.  Equivalently,  $\ell\in \widetilde{I}_1$ if  and only if   $\xi_{K_\ell} >0$, where $K_ \ell=  \llbracket 1,p\rrbracket-\{\ell\}  $. Hence  $\widetilde{I}_1$ can be explicited by sorting the  $\xi_{K_\ell} $s in decreasing order since this approach makes it possible to identify the non-zero terms. Indeed, since  $\widetilde{I}_0$ is not empty,   there exist integers $\nu_1,\cdots,\nu_p$ and   $\mathcal{D}\in\{1,\cdots,p-1\}$ such that:
\begin{equation}\label{tri}
\xi_{K_{\nu_1}}\geq\xi_{K_{\nu_2}}\geq\cdots\geq\xi_{K_{\nu_\mathcal{D}}}>0=\xi_{K_{\nu_{\mathcal{D}+1}}}=\cdots=\xi_{K_{\nu_{p}} },
\end{equation}
with $\nu_i<\nu_\ell$ if $\xi_{K_i}=\xi_{K_\ell}$ and $i<\ell$. Hence
  $\widetilde{I}_1$ can be writen as

\begin{equation}\label{i1}
\widetilde{I}_1=\{\nu_1,\cdots,\nu_\mathcal{D}\}.
\end{equation}
and its estimation  reduces to that of the parameters  $\nu_\ell$  and $\mathcal{D}$. Our method for selecting variables will be based on estimating these parameters.

\section{Selection of variables}\label{sec:MLE}
In this section we propose estimates  of the aforementioned parameters by using an approach tackled in Mbina Mbina et al. (2023)  in the context of multivariate linear regression, so achieving variable selection via the resulting estimation of   $\widetilde{I}_1$.   This later estimate  depends of tuning parameters, so a procedure for choosing optimal values for these parameters, based on $V$-fold cross validation, is introduced.

\subsection{Estimation}\label{sec:MLE1}
From the sample  $\left\{\left(Y^{(i)},X^{(i)}\right)\right\}_{1\leq i\leq n}$ we consider the samples $\{\textrm{\textbf{X}}_{\centerdot\ell}^{(i)}\}_{1\leq i\leq n}$  of coordinates of the $X_{\ell}^{(i)}$'s on the basis $\{\phi_{k\ell}\}_{1\leq k\leq d_\ell}$, that is 
\begin{equation*}\label{coord}
X_{\ell}^{(i)}(t) =\left( \textrm{\textbf{X}}_{\centerdot\ell}^{(i)}\right)^T\phi_{\centerdot\ell}(t),
\end{equation*}
and we put
\begin{equation}\label{Xi}
\mathcal{X}^{(i)}= \left(
 \begin{array}{c}
\textrm{\textbf{G}}_{1}\, \textrm{\textbf{X}}^{(i)}_{\centerdot 1}\\
\hline 
 \textrm{\textbf{G}}_{2}\, \textrm{\textbf{X}}^{(i)}_{\centerdot 2} \\
\hline 
 \vdots \\
\hline 
\textrm{\textbf{G}}_{p}\, \textrm{\textbf{X}}^{(i)}_{\centerdot p}
 \end{array}
\right) .
\end{equation}
We then consider the sample means
\begin{equation}\label{means}
\overline{\mathcal{X}}^{(n)} = \frac{1}{n}\sum_{i=1}^{n}\mathcal{X}^{(i)}, \hspace{0.2cm} \overline{Y}^{(n)} =  \frac{1}{n}\sum_{i=1}^{n}Y^{(i)},
\end{equation}
and the empirical covariance and cross-covariance  matrices
\begin{equation}\label{c1}
\widehat{ \mathcal{C}}_1 =  \frac{1}{n}\sum_{i=1}^{n}(\mathcal{X}^{(i)}- \overline{\mathcal{X}}^{(n)})\,(\mathcal{X}^{(i)}- \overline{\mathcal{X}}^{(n)})^T, 
\end{equation}
and
\begin{equation}\label{c12}
\widehat{ \mathcal{C}}_{12}=  \frac{1}{n}\sum_{i=1}^{n}(\mathcal{X}^{(i)} - \overline{\mathcal{X}}^{(n)})(Y^{(i)} - \overline{Y}^{(n)})^T
\end{equation}
from which we  estimate  $\xi_K$  by

\begin{equation*}\label{crit2} 
\widehat{\xi}_{K}= \| \widehat{ \mathcal{C}}_{12}- \widehat{ \mathcal{C}}_1\widehat{\Pi}_K \widehat{ \mathcal{C}}_{12}\|
\end{equation*}
where $$\widehat{\Pi}_K = A_{K}^{T}(A_{K}\widehat{ \mathcal{C}}_1A_{K}^{T})^{-1}A_{K}.$$
A naive approach for estimating the $\nu_\ell$s consists in  sorting the $\widehat{\xi}_{K_\ell}$s in decreasing order as it was done in \eqref{tri} with the   $\xi_{K_\ell}$s, but since  such an approach does not  guarantee the consistency of the resulting estimator  because  of possible ties,  we will rather use an  estimate of  $\xi_{K_\ell}$  obtained from an appropriate penalization  of  $\widehat{\xi}_{K_\ell}$  which allows to avoid ties, so yielding consistency.   More specifically, we consider the statistics
\begin{equation*}\label{pen1}
\widehat{\phi}_{\ell}=\widehat{\xi}_{K_{\ell}}+\frac{f\left(  \ell\right)}{n^{\alpha}}, 
\end{equation*}
 where
$0<\alpha<1/2  $ and $f$ is  a strictly decreasing function from
$ \llbracket 1,p\rrbracket$ to $\mathbb{R}_{+}$,  that we sort in decreasing order so as to obtain the integers  $\widehat{\nu}_1\cdots, \widehat{\nu}_p$ belonging to $\llbracket 1,p\rrbracket$,  satisfying
\begin{equation*}\label{decroissant}
\widehat{\phi}_{\widehat{\nu}_1}>\widehat{\phi}_{\widehat{\nu}_2}>\cdots>\widehat{\phi}_{\widehat{\nu}_p}
\end{equation*}
and which estimate $\nu_1\cdots, \nu_p$ respectively.  In order to estimate $\mathcal{D}$ we will first highlight a property that characterizes this parameter. Considering   the subset  $J_\ell=\{\nu_k\,/\,1\leqslant k\leqslant \ell\}$ of $\llbracket 1,p\rrbracket$, it is seen  from \eqref{i1} that $J_\ell\subset\widetilde{ I}_1$ if  $\ell\geqslant \mathcal{D}$. This imples that    $\xi_{ J_\ell}=0$ if $\ell\geqslant \mathcal{D}$, and  $\xi_{ J_\ell}>0$ if $i< \mathcal{D}$ and shows, therefore, that $\mathcal{D}$ is the   smallest integer $\ell\in \llbracket 1,p\rrbracket$ for which $\xi_{ J_\ell}$ has its minimum value. So, for estimating $\mathcal{D}$ we will  minimize   an estimate of  $\xi_{ J_\ell}$. More precisely, for the same reason than above, we will use a penalized estimate of this index  obtained as  
\begin{equation*}\label{pen2}
\widehat{\psi}_{\ell}=\widehat{\xi}_{\widehat{J}_{\ell}}+\frac{g\left( \widehat{\nu}_\ell\right)}{n^{\beta}}, 
\end{equation*}
where $\widehat{J}_{\ell}=\left\{
\widehat{\nu}_1,\cdots,\widehat{\nu}_\ell\right\}  $, $0<\beta<1/2  $ and $g$ is  a strictly increasing function  from $ \llbracket 1,p\rrbracket$ to \ $\mathbb{R}_{+}$.
Then, we estimate $\mathcal{D}$ by
\[
\widehat{\mathcal{D}}=\arg\min_{\ell\in \llbracket 1,p\rrbracket}\left(  \widehat{\psi}_{\ell}\right)   
\]
and take  the set
\[
\widehat{I}_{1}=\left\{ \widehat{\nu}_1,\widehat{\nu}_2,\cdots,\widehat{\nu}_{ \widehat{\mathcal{D}}}\right\}
\]
as the required set of indices of the relevant explanatory functional variables in model \eqref{model}.

\subsection{Choosing optimal tuning parameters}\label{sec:MLE2}
The procedure for variable selection introduced in the preceding section depends on two tuning parameters  $\alpha$ and $\beta$ which   may have influence on the performance of our method; then choosing optimal values for  these parameters  is a crucial issue.  We propose an optimal choice of  $(\alpha,\beta)$  based  on $V$-fold  cross validation (with $V\in\mathbb{N}^\ast$) used in order to minimize the mean squared error of prediction (MSEP), that is  a distance between the $Y^{(i)}$'s and their predictions $\widehat{Y}^{(i)}$ obtained by least squared method from model \eqref{modelfinal}.  Since this model reduces to
\[
Y=\text{\textbf{B}}^T\mathcal{X}+\varepsilon,
\]
where $\text{\textbf{B}}$ is the $d\times q$ matrix given by
\[
\text{\textbf{B}}=\left(
 \begin{array}{c}
b_{\centerdot \centerdot 1}\\
\hline 
b_{\centerdot \centerdot 2}\\
\hline 
 \vdots \\
\hline 
b_{\centerdot \centerdot p}
 \end{array}
\right) ,
\]
we have  $\widehat{Y}^{(i)}= \widehat{\text{\textbf{B}}}^T\mathcal{X}^{(i)}$, where
\[
 \widehat{\text{\textbf{B}}}=\bigg(\sum_{i=1}^n\mathcal{X}^{(i)} \mathcal{X}^{(i)^T}\bigg)^{-1}\bigg(\sum_{i=1}^n\mathcal{X}^{(i)} Y^{(i)^T}\bigg).
\]
Considering the $n\times d$ and $n\times q$ matrices $\mathbb{X}$ and $\mathbb{Y}$ defined as
\[
\mathbb{X}^T=\left(
 \begin{array}{ccccc}
\mathcal{X}^{(1)} & \big\vert & \cdots & \big\vert & \mathcal{X}^{(n)}
 \end{array}
\right)\,\,\,\textrm{ and }\,\,\,\mathbb{Y}^T=\left(
 \begin{array}{ccccc}
Y^{(1)} & \big\vert & \cdots & \big\vert & Y^{(n)}
 \end{array}
\right),
\]
we have $ \widehat{\text{\textbf{B}}}=\left(\mathbb{X}^T\mathbb{X}\right)^{-1}\mathbb{X}^T\mathbb{Y}$, and the MSEP  is 
\[
\textrm{MSEP}=\frac{1}{n}\sum_{i=1}^n\Vert   Y^{(i)}-\widehat{Y}^{(i)}\Vert_{\mathbb{R}^q}^2=\frac{1}{n}\Vert \mathbb{Y}-\mathbb{X}\left(\mathbb{X}^T\mathbb{X}\right)^{-1}\mathbb{X}^T\mathbb{Y}\Vert^2.
\]
When a subset $K=\{i_1,\cdots,i_k\}$ of explanatory variables is used for the prediction, then  the MSEP  computed on a subsample  $\left\{\left(Y^{(i)},X^{(i)}\right)\right\}_{  i\in S}$ , where $S$ is  the subset  $S=\{s_1,\cdots,s_m\}$  of indices in $ \llbracket 1,n\rrbracket$,  is
\begin{equation}\label{cval}
\textrm{MSEP}_K^S=\frac{1}{m}\Vert \mathbb{Y}_S-\mathbb{X}_SA_K^T\left(A_K\mathbb{X}_S^T\mathbb{X}_SA_K^T\right)^{-1}A_K\mathbb{X}_S^T\mathbb{Y}_S\Vert^2,
\end{equation}
where
\[
\mathbb{X}^T_S=\left(
 \begin{array}{ccccc}
\mathcal{X}^{(s_1)} & \big\vert & \cdots & \big\vert & \mathcal{X}^{(s_m)}
 \end{array}
\right) \,\,\,\textrm{ and }\,\,\,\mathbb{Y}^T_S=\left(
 \begin{array}{ccccc}
Y^{(s_1)} & \big\vert & \cdots & \big\vert & Y^{(s_m)}
 \end{array}
\right).
\]
Now, consider a partition $\{\mathcal{S}_1,\cdots , \mathcal{S}_V\}$ of the set  $\mathcal{S}=\{1,\cdots,n\}$, each $\mathcal{S}_j$ having the same size $m\in\mathbb{N}^\ast$ (then,   $n=mV$). For each $j$ in $\{1,\cdots,V\}$, after removing the $j$-th subset  $\mathcal{S}_j$ from $\mathcal{S}$,  we apply our method for selecting variable on the remaining  subsample  $\left\{\left(Y^{(i)},X^{(i)}\right)\right\}_{  i\in \mathcal{S}-\mathcal{S}_j}$   with a given value for $(\alpha,\beta)$ in $]0,1/2[^2$; this leads to  an  estimate  $\widehat{I}_1^{(-j)}$ of $\widetilde{I}_1$. Then, we  define the  cross-validation index
\[
CV(\alpha,\beta)=\frac{1}{V}\sum_{j=1}^V\textrm{PL}^{(j)}(\alpha,\beta),
\]
where $\textrm{PL}^{(j)}(\alpha,\beta)$ is the MSEP given in \eqref{cval} with $K=\widehat{I}_1^{(-j)}$ and $S=\mathcal{S}_j$. An optimal value  $(\widehat{\alpha},\widehat{\beta})$   of   $(\alpha,\beta)$  is obtained by minimizing this index, that is
\begin{equation}\label{opti}
(\widehat{\alpha},\widehat{\beta})=\underset{(\alpha,\beta)\in ]0,1/2[^2}{\mathrm{argmin}}CV(\alpha,\beta).
\end{equation}
\subsection{Algorithms}
In this section we give  algorithms from which our proposal can concretely be computed. Three algorithms are presented. The first one describes the proposed approach for determining the dimensions for basis representations. The second algorithm gives the steps for obtaining $\widehat{I}_1$ from a given sample and a given value of the pair $(\alpha,\beta)$. The third algorithm shows how to use the previous ones to implement our variable selection method.

\bigskip

\begin{algorithm}\label{calcdim}
\caption{Computation of $d_\ell$, $\ell=1,\cdots,p$.}\label{algo1}
\hspace*{\algorithmicindent} \textbf{Input:} a sample $\{X^{(i)}\}_{1\leqslant i\leqslant n}$ of functional variables as in \eqref{sample}; fine grid of points  

\hspace{1.8cm}$t^{(\ell)}_1,\cdots,t_{N_\ell}^{(\ell)}$  in $\mathcal{I}_ \ell$, $\ell=1,\cdots,p$; bases $\mathscr{B}_\ell=\{\phi_{k\ell}\}_{k\geqslant 1}$ of  $L^2(\mathcal{I}_\ell)$, $\ell=1,\cdots,p$;

\hspace{1.8cm}a maximal dimension $d_\textrm{max}$.\\
\hspace*{\algorithmicindent} \textbf{Output:} optimal dimensions $d_1\cdots,d_p$.
\medskip
\begin{algorithmic}[1]
\FOR {$\ell=1,\cdots, p$}
\FOR {$i=1,\cdots, n$}
\FOR {$m=1,\cdots, d_\textrm{max}$}
\STATE  compute the coordinate of $\mathbf{X}^{(i)}_{k\ell}$ of $X^{(i)}_\ell$ on $\phi_{k\ell}$, for $k=1,\cdots,m$
\STATE compute $\textrm{BIC}_m(i,\ell)=\ln\left(\sum_{r=1}^N\left(X^{(i)}_\ell(t^{(\ell)}_r)- \sum_{k=1}^{m} \textrm{\textbf{X}}^{(i)}_{k \ell}\phi_{k \ell}(t_r^{(\ell)})\right)^2\right)$ 

\hspace{4cm}$+(m+1)\left(\frac{\ln (N_\ell)}{N_\ell}\right)$
\ENDFOR
\STATE set $d^{(i)}_\ell=\arg\min_{1\leqslant m\leqslant d_\textrm{max}}\left(\textrm{BIC}_m(i,\ell)\right)$
\ENDFOR
\STATE set $d_\ell=\max\{d^{(1)}_\ell,\cdots,d^{(p)}_\ell\}$
\ENDFOR
\end{algorithmic}
\end{algorithm}

\begin{algorithm}\label{calci1}
\caption{Computation of  $\widehat{I}_1$}\label{algo2}
\hspace*{\algorithmicindent} \textbf{Input:} a sample $\{(X^{(i)},Y^{(i)})\}_{1\leqslant i\leqslant n}$  as in \eqref{sample}; fine grid of points  

\hspace{1.8cm}$t^{(\ell)}_1,\cdots,t_{N_\ell}^{(\ell)}$  in $\mathcal{I}_ \ell$, $\ell=1,\cdots,p$; bases $\mathscr{B}_\ell=\{\phi_{k\ell}\}_{k\geqslant 1}$ of  $L^2(\mathcal{I}_\ell)$, $\ell=1,\cdots,p$;

\hspace{1.8cm}a pair $(\alpha,\beta)$ of tuning parameters belonging to $]0,1/2[^2$; penalty functions 

\hspace{1.8cm}$f$ and $g$; Gram matrices  $\mathbf{G}_1,\cdots,\mathbf{G}_p$.\\
\hspace*{\algorithmicindent} \textbf{Output:} subset $\widehat{I}_1$ of    selected variables indices.
\medskip
\begin{algorithmic}[1]
\STATE compute the dimensions $d_1,\cdots,d_p$ by using Algorithm \ref{algo1} on the sample  $\{X^{(i)}\}_{1\leqslant i\leqslant n}$
\FOR {$i=1,\cdots, n$}
\FOR {$\ell=1,\cdots,p$}
\STATE  compute the coordinate $\mathbf{X}^{(i)}_{\centerdot\ell}$ of $X^{(i)}_\ell$ on the basis  $\mathscr{B}_\ell$
\STATE  compute the Gram matrix $\mathbf{G}_\ell$ according to \eqref{gfourier}, \eqref{gspline} or \eqref{ggaussian}; in case of Fourier basis set  $\mathbf{G}_\ell=\mathbb{I}_{d_\ell}$
\ENDFOR
\STATE  compute   $\mathcal{X}^{(i)}$ as defined in \eqref{Xi}
\ENDFOR
\STATE  compute the sample means, covariance and cross-covariance matrix given in \eqref{means}, \eqref{c1} and \eqref{c12}
\FOR  {$\ell=1,\cdots,p$}
\STATE compute $\widehat{\phi}_\ell=\widehat{\xi}_{K_{\ell}}+n^{-\alpha} f\left(  \ell\right)$
\ENDFOR
\STATE set $\widehat{\nu}_1,\cdots,\widehat{\nu_p}$ that  satisfy $\widehat{\phi}_{\widehat{\nu}_1}>\widehat{\phi}_{\widehat{\nu}_2}>\cdots>\widehat{\phi}_{\widehat{\nu}_p}$
\FOR  {$\ell=1,\cdots,p$}
\STATE set $\widehat{J}_\ell=\{\widehat{\nu}_1,\cdots,\widehat{\nu_\ell}\}$ 
\STATE compute $\widehat{\psi}_{\ell}=\widehat{\xi}_{\widehat{J}_{\ell}}+n^{-\beta} g\left( \widehat{\nu}_\ell\right)$ 
\ENDFOR
\STATE set  $\widehat{\mathcal{D}}=\arg\min_{\ell\in \llbracket 1,p\rrbracket}\left(  \widehat{\psi}_{\ell}\right)  $
\STATE   set  $\widehat{I}_{1}=\left\{ \widehat{\nu}_1,\widehat{\nu}_2,\cdots,\widehat{\nu}_{ \widehat{\mathcal{D}}}\right\}$ 
\end{algorithmic}
\end{algorithm}

\begin{algorithm}\label{selectionne}
\caption{The proposed method for variable selection}\label{algo3}
\hspace*{\algorithmicindent} \textbf{Input:} a sample $\{(X^{(i)},Y^{(i)})\}_{1\leqslant i\leqslant n}$  as in \eqref{sample}; fine grid of points  

\hspace{1.8cm}$t^{(\ell)}_1,\cdots,t_{N_\ell}^{(\ell)}$  in $\mathcal{I}_ \ell$, $\ell=1,\cdots,p$; bases $\mathscr{B}_\ell=\{\phi_{k\ell}\}_{k\geqslant 1}$ of  $L^2(\mathcal{I}_\ell)$, $\ell=1,\cdots,p$;

\hspace{1.8cm}penalty functions $f$ and $g$.\\
\hspace*{\algorithmicindent} \textbf{Output:} subset $\widehat{I}_1$ of    selected variables indices.
\begin{algorithmic}[1]
\STATE  divide  the sample  $\{(X^{(i)},Y^{(i)})\}_{1\leq i\leq n}$ into two subsamples : a training sample $S_\mathscr{L}$ and a test sample $S_\mathscr{T}$
\STATE  partition  $S_\mathscr{L}$ into $V$ subsamples  $ \mathcal{S}_1,\cdots , \mathcal{S}_V $   having the   size $m $ such that $n=mV$ 
\STATE  get a  grid $\mathcal{G}$ of pairs  from a fine discretization of   $ ]0,1/2[^2$  
\FOR {$(\alpha,\beta)\in \mathcal{G}$}
\FOR  {$j=1,\cdots,V$}
\STATE remove $ \mathcal{S}_j$  from $S_\mathscr{L}$, the remaining data set is denoted by $\mathcal{S}_\mathscr{L}^{(-j)}$
\STATE  get  $\widehat{I}_1^{(-j)}$ by applying  the   method of Algorithm \ref{algo2} on  $\mathcal{S}_\mathscr{L}^{(-j)}$ with $(\alpha,\beta)$
\STATE  compute $\textrm{PL}^{(j)}(\alpha,\beta)$ by using  \eqref{cval} with $K=\widehat{I}_1^{(-j)}$ and $S=\mathcal{S}_j$
\ENDFOR
\STATE  set  $CV(\alpha,\beta)=\frac{1}{V}\sum_{j=1}^V\textrm{PL}^{(j)}(\alpha,\beta)$
\ENDFOR
\STATE  take $(\widehat{\alpha},\widehat{\beta})$ that minimizes $CV(\alpha,\beta)$ over $\mathcal{G}$
\STATE   get  $\widehat{I}_1$ by applying  the   method of Algorithm \ref{algo2}  on  $\mathcal{S}_\mathscr{T}$  with $(\widehat{\alpha},\widehat{\beta})$ 
\end{algorithmic}
\end{algorithm}

\section{Simulations}\label{sec4:sim}
\noindent We investigate  through simulations the   performance  of the proposed variable selection method that we compare to those of the random subspace method   of Smaga and Matsui (2018)   and the group SCAD method of  Matsui and Konishi (2011).  The  data sets  was generated according to the following examples:

\begin{itemize}
  \item \textbf{Example 1:} We set   $ p=10$,  $q=1$ and we generate independently  the processes $X_\ell(t)= 5\sum_{k=1}^{50} c_{\ell}^{(k)}\psi_{\ell}^{(k)}(t)$, $t\in [0,1]$, $\ell=1,\cdots,10$, where  $c_{\ell}^{(k)} \sim N(0,k^{-2})$, $\psi_{\ell}^{(1)}(t)= 1$ and $\psi_{\ell}^{(k+1)}(t)$ = $\sqrt{2}\cos(k\pi t)$ for $k\geq 1$.   The functional coefficients are taken as $B_\ell(t)=b_\ell\sin(\pi\ell t/10)$   with $b_1= 0.25$, $b_5= 0.50$,  $b_6= 0.75$,  $b_{7}= 1.00$, $b_{10}=  1.25$, and $b_2=b_3=b_4=b_8=b_9=0$. The  error term is taken as  $\varepsilon \sim$ N(0, $\sigma^{2}$) with $\sigma=0.1$, $0.25$, $0.5$. The response $Y$ is then generated  according to model  \eqref{model} with the related  integrals    computed by using the trapezoidal rule on the basis of the values of the $X_\ell$s and the $B_\ell$s on $N=51$ equidistant points $t_j=j/50$, $j=0,\cdots,50$, in $[0,1]$.  In this example, the true set of relevant variables is $I_{1}=\left\{1,5,6,7,10\right\}$.
\end{itemize}

\begin{itemize}
  \item \textbf{Example 2:} We set   $ p=6$,  $q=1$ and we generate independently  the processes defined for $t\in [0,1]$ by:

\smallskip

$X_1(t)= a_1t^3+a_2t^2+a_3t+a_4$ with  $a_1\sim N(-2,1)$, $a_2\sim U(2,3)$, $a_3\sim Exp(1)$, $a_4\sim N(0,0.1)$,

\smallskip

$X_2(t)=b_1\sin(2\pi t/3)+b_2t$ with $b_1\sim U(3,7)$, $b_2\sim N(0,1)$,

\smallskip
 
$X_3(t)=c_1(2t-1)^3+c_2(2t-1)^2+c_3(2t-1)+c_4$ with  $c_1\sim N(-3,1.2)$,
 $c_2\sim N(2,0.5)$, $c_3\sim N(-2,1)$, $c_4\sim N(2,1.5)$,

\smallskip
 
$
X_4(t)=(t-d_1)^2\cos(2\pi t/3)+d_2t+d_3$ with $d_1\sim U(2,1)$, $d_2\sim N(0,1)$, $d_3\sim Exp(1)$,

\smallskip

$X_5(t)=\cos(2\pi (t-e_1))+e_2t+e_3$ with $e_1\sim N(-5,3)$, $e_2\sim N(7,1)$, $e_3\sim N(0,0.025)$,

\smallskip
and

\smallskip

$X_6(t)=f_1t^8+ \cos(f_2\pi t)+ t^4\sin(f_3\pi t)+f_4$ with  $f_1\sim N(-4,2)$, $f_2\sim U(0,1)$,
$f_3\sim U(0,1/2)$, $f_4\sim N(0,0.1)$.

\bigskip

The functional coefficients are taken as   $B_1(t)= t\sin(\pi t/4)$,  $B_2(t)= \cos(2\pi t)+t^2+1$,  $B_5(t)=e^{-2t}+t^3-1$, and $B_3(t)=B_4(t)=B_6(t)=0$. Here also,  we take the error term as   $\varepsilon \sim$ N(0, $\sigma^{2}$) with $\sigma=0.1$, $0.25$, $0.5$, and  the response $Y$ is  generated as in Example 1.  In this example, the true set of relevant variables is $I_{1}=\left\{1,2,5\right\}$. 
\end{itemize}

\begin{itemize}  
	\item \textbf{Example 3:}  We set   $ p=8$,  $q=2$ and we generate independently  the processes defined  on $N=51$ equidistant points, $t_s=j/50$, $s=0,\cdots,50$, in $[0,1]$ by $X_\ell(t_s)=u_\ell(t_s)+\eta_{\ell s}$, $\ell=1,\cdots,8$, where  

$u_{1}(t)=a_{1}(2t-1)^3+a_{2}(2t-1)^2+a_{3}(2t-1)+a_{4}$ with $a_{1}\sim \mathcal{N}(-3,1.2)$, $a_{2} \sim \mathcal{N}(2,0.5),$ $a_{3}\sim \mathcal{N}(-2,1)$, $a_{4} \sim \mathcal{N}(2,1.5)$,

\smallskip

$u_{2}(t)=b_{1}t^8+\cos(b_{2}\pi t)+b_{3}t^4\sin(b_3\pi t)+b_{4}$ with $b_{1}\sim N(-4,2)$, $b_{2} \sim U(0,1),$ $b_{3}\sim U[0,0.5]$, $b_{4} \sim N(0,0.1),$

\smallskip

$u_{3}(t)=c_{1}\cos(2\pi t)+c_{2}$ with  $c_{1} \sim N(-4,3)$, $c_{2} \sim N(7,1.5),$

\smallskip

$u_{4}(t)=d_{1}\sin(\pi^2 t/3)+d_{2}$ with  $d_{1}\sim U(3,7)$, $d_{2} \sim$ $N(0,1),$

\smallskip

$u_{5}(t)=e_{1}\cos^3(3\pi (2t-1)+e_{2}\cos^2(2\pi (2t-1))+e_{3}\cos^3(\pi (2t-1))$ with  $e_{1} \sim N(-3,1.2)$, $e_{2} \sim N(2,0.5)$, $e_{3} \sim N(-2,1),$

\smallskip
$u_{6}(t=f_{1}\sin(2\pi^2 t/3)+f_{2}\cos(\pi^2 t/3)$ with  $f_{1}$ $\sim$ $N(-2,1)$, $f_{2}$ $\sim N(3,1.5);$

\smallskip

$u_{7}(t)=g_{1}\cos(2\pi (3t-2))$ + $g_{2}\cos(\pi (3t-2))$ with  $g_{1}\sim U(2,7)$, $g_{2} \sim N(2,0.4),$

\smallskip

$u_{8}(t)=h_{1}\cos(\pi (2t-1))+h_{2}(2t-1)+h_{3}$ with  $h_{1} \sim N(4,2)$, $h_{2} \sim N(-3,0.5)$, $h_{3} \sim N(1,1),$

and  $\eta_{\ell s}\sim N(0, 0.025\nu^2_{\ell s})$  with $\nu^2_{\ell s}=\max_{0\leqslant s\leqslant 50}\{u_{\ell}(t_{s})\} - \min_{0\leqslant s\leqslant 50}\{u_{\ell}(t_{s})\}$. The functional coefficients are taken as   $B_{13}(t)= 0.25\sin(t)$,  $B_{15}(t)=  0.75\sin(2t-1)$, $B_{17}(t)=  1.25\sin(3t-2)$,   $B_{23}(t)= 0.25\cos(t)$,  $B_{25}(t)=  0.75\cos(2t-1)+(2t-1)^2$, $B_{27}(t)=  1.25\cos(3t-2)+(3t-2)^4$ and $B_{j\ell}(t)=0$ for $j=1,2$, $\ell=1,2,4,6,8$. The  errors  are generated independently as    $\varepsilon_j \sim N(0, \sigma^{2})$, $j=1,2$,    $\sigma=0.1$, $0.25$, $0.5$, and the reponses $Y_1$ and $Y_2$ are generated    according to model  \eqref{model} with   integrals    computed by using the trapezoidal rule based on  the  equidistant points $t_s$ introduced above.  In this example, the true set of relevant variables is $I_{1}=\left\{3,5,7\right\}$.
\end{itemize}

\bigskip

\noindent We simulate $200$ independent replications of samples from the models given    in the above examples. For each replication:
\begin{itemize}
\item a training sample of size $n=50$, $75$, $100$ is generated, and is used for computing optimal values of the tuning parameters on which our method lies by using $V$-fold cross validation, with $V=5$,  as described in Section \ref{sec:MLE2};
\item  a test  sample having the same  size   is generated. On this sample, our method for variable selection  is performed with the aforementioned optimal values of tuning parameters, together with the random subspace method   of Smaga and Matsui (2018) and the group SCAD method of   Matsui and Konishi (2011).
\end{itemize}
Over these $200$ replications, the following four measurements are computed in order to assess performance of the different methods:
\begin{itemize}
\item[$(1)$]average model size, i.e.  MSIZE$=200^{-1}\sum_k\vert \widehat{I}^{(k)}\vert$, where  $\widehat{I}^{(k)}$ is the subset of selected variables at the $k$th replication;
\item[$(2)$]coverage probability, i.e.  CVP$=200^{-1}\sum_k I(I_1\subset  \widehat{I}^{(k)})$, where  $I(\cdot)$ denotes the indicator function;
\item[$(3)$]average false discovery rate, i.e. FDR$=200^{-1}\sum_kN_k/\vert \widehat{I}^{(k)}\vert$, where  $N_{k}$ is the  number of false discovery variables for the $ k$th replication;
\item[$(4)$]mean squared errors  of prediction (MSEP) computed on the test sample  after variable selection , i.e. $\textrm{MSEP}_K^S$ given  in  \eqref{cval}, with  $S= \llbracket 1,n\rrbracket$ and $K= \widehat{I}^{(k)}$, for $k=1,\cdots,200$.  
\end{itemize}

\bigskip

\noindent  We use  the R programming language for performing the three  methods. The basis representations of the  functional explanatory variables  are estimated by using  the following  R functions of the package \textit{fda}: \textit{create.fourier.basis()},  \textit{create.spline.basis()}, \textit{eval.basis()}, \textit{fdPar()} and \textit{smooth.basis()}. The   R the  function \textit{grpreg()}  of  the package \textit{grpreg} is used  for performing the group SCAD method. 

\bigskip

\noindent Tables \ref{tab1}, \ref{tab2} and \ref{tab3}  report the obtained results for CVP, FDR and MSIZE, whereas those  related to MSEP are given in Figures 1  to  6. In Table \ref{tab1},  the group SCAD method  (denoted by gSCAD)   outperforms the two others  with regard to CVP, but our method gives  better  results than it  concerning FDR and MSIZE.  On the other hand,  Table 2 shows a slight superiority of our method over gSCAD  with regard to CVP, these two methods being much better than the random subspace  method (denoted by RSM). Nevertheless, it can be observed  in this table that  gSCAD gives better results than ours and  RSM   in term of  FDR. Table 3 shows the performance of our method in the multivariate case ($q>1$), and it is alone there because the other two methods are not appropriate for this case. We observe good results with sufficiently high values of CVP and   low values of FDR. Concerning MSEP, we can see very low values for it in almost all cases, with a superiority of our method over the two others in Figure 3 and Figure 4. Comparing the basis representations, we see in  Figure 5 and Figure 6  that better results are obtained with the B-spline  basis  for low sample size ($n=50$), and by the Fourier basis for large sample size ($n=100$).

\bigskip

\begin{center}
\begin{table}
\caption{Results of our method (OM), random subspace method (RSM) and group SCAD method (gSCAD)  in term of  the  coverage probability (CVP), false discovery rate (FDR) and model size (MSIZE) across $200$ replications for Example 1.}\label{tab1}
\medskip
\begin{tabular}{cccccccccccccccccc}%
\hline\hline
    &  & &        &        &   &  &\multicolumn{3}{c }{Fourier basis}&     &  & &     & \multicolumn{3}{c }{B-spline basis}&  \\
\cline{8-10}\cline{15-17}
  & n & & & $\sigma$ & Method   &   & CVP  & FDR & MSIZE  &  &   &  & &     CVP & FDR & MSIZE  & \\
\hline
& 50 & & &0.10 & OM     &  & 0.21 & 0.00 & 4.00 &   & &  &  & 0.26 & 0.00 & 4.00 &  \\
 &   & & &    & RSM   &  & 0.25 & 0.22 & 5.50 &  &  &  & &  0.50 & 0.37 & 5.50  & \\
&    &  & &   & gSCAD  &   & 0.75 & 0.51 & 9.75 &  &  &  & &  0.25 & 0.50 & 10.0 &  \\
 & & &    &0.50 & OM   &  & 0.37 & 0.00 & 4.00 &   &   &  & &  0.36 & 0.00 & 4.00 &\\
 &  & & &   & RSM   &   & 0.22 & 0.48 & 6.75 &  &  &   & & 0.37 & 0.52 & 6.50 & \\
  &  &  & &   & gSCAD   &  & 1.00 & 0.53 & 10.0 &   &   & &  &1.00 & 0.54 & 10.0 & \\
\hline
& 75 & & &0.10 & OM     &  & 0.42 & 0.38 & 4.60 &   & &  &  & 0.20 & 0.25 & 4.00 &  \\
 &   & & &    & RSM   &  & 0.50 & 0.20 & 5.75 &  &  &  & &  0.17 & 0.12 & 3.25  & \\
&    &  & &   & gSCAD  &   &1.00 & 0.50 & 9.75 &  &  &  & &  1.00 & 0.50 & 10.0 &  \\
 & & &    &0.50 & OM   &  & 0.33 & 0.25 & 3.75 &   &   &  & &  0.46 & 0.37  & 3.58 &\\
 &  & & &   & RSM   &   & 0.25 & 0.43 & 6.50 &  &  &   & & 0.50 & 0.34 & 5.50 & \\
  &  &  & &   & gSCAD   &  & 1.00 & 0.50 & 10.0 &   &   & &  &1.00 & 0.50 & 10.0 & \\
\hline
 &100 & & &0.10 & OM   &   & 0.31 & 0.43 & 4.70  &   &   & &  & 0.34& 0.12 & 3.84 &\\
  &   & & &    & RSM   &  & 0.25 & 0.30 & 6.50  &  &  &   & & 0.20 & 0.28 & 4.50 & \\
  &  & & &    & gSCAD   &   & 1.00 & 0.50 & 9.85 &   &   &   & & 1.00 & 0.50 & 10.0 &\\
  &  &  &  &0.50 & OM   &   & 0.39 & 0.00 & 4.50  &  &   &  & & 0.41 & 0.11 &4.80 & \\
    &  &  &   &    & RSM   &  & 0.25 & 0.42 &5.75  &  &  &  & & 0.31 & 0.17 & 6.35  & \\
   &  &  & &    & gSCAD    &   & 1.00 & 0.50 & 9.90 &   &   & & &1.00 & 0.50 & 9.95 & \\
\hline\hline
\end{tabular}
%\end{minipage}
\end{table}
\end{center}

\bigskip

\begin{center}
\begin{table}
\caption{Results of our method (OM), random subspace method (RSM) and group SCAD method (gSCAD)  in term of  the  coverage probability (CVP), false discovery rate (FDR) and model size (MSIZE) across $200$ replications for Example 2.}\label{tab2}
\medskip
\begin{tabular}{cccccccccccccccccc}%
\hline\hline
    &  & &        &        &   &  &\multicolumn{3}{c }{Fourier basis}&     &  & &     & \multicolumn{3}{c }{B-spline basis}&  \\
\cline{8-10}\cline{15-17}
  & n & & & $\sigma$ & Method   &   & CVP  & FDR & MSIZE  &  &   &  & &     CVP & FDR & MSIZE  & \\
\hline
& 50 & & &0.10 & OM     &  & 1.000 & 0.455 & 5.600 &   & &  &  & 1.000 & 0.455 & 5.500 &  \\
 &   & & &    & RSM   &  & 0.000 & 0.000 & 1.000 &  &  &  & &  0.000 & 0.000 & 1.000  & \\
&    &  & &   & gSCAD  &   & 1.000 & 0.013 & 3.050 &  &  &  & &  0.900 & 0.130 & 2.950 &  \\
 & & &    &0.50 & OM   &  & 0.850 & 0.405 & 5.200 &   &   &  & & 0.850 & 0.420 &5.350 &\\
 &  & & &   & RSM   &   & 0.000 & 0.000 &1.000 &  &  &   & & 0.000 & 0.000 & 1.000 & \\
  &  &  & &   & gSCAD   &  & 0.800 & 0.367 & 4.450 &   &   & &  &0.700 & 0.140 & 3.200 & \\
\hline
& 75 & & &0.10 & OM     &  & 1.000 & 0.475 & 5.750 &   & &  &  &1.000 & 0.455 & 5.550 &  \\
 &   & & &    & RSM   &  & 0.000 & 0.000 & 1.119 &  &  &  & &  0.000 & 1.000 & 1.000  & \\
&    &  & &   & gSCAD  &   & 1.000 & 0.140 & 3.700 &  &  &  & &  1.000 & 0.000 & 3.000 &  \\
 & & &    &0.50 & OM   &  & 1.000 & 0.471 & 5.650 &   &   &  & &  0.950 & 0.442 & 5.500 &\\
 &  & & &   & RSM   &   & 0.000 & 0.000 &1.000 &  &  &   & & 0.000 & 0.000 & 1.000 & \\
  &  &  & &   & gSCAD   &  & 0.950 & 0.329 & 4.550 &   &   & &  &0.900 & 0.195 & 3.850 & \\
\hline
& 100 & & &0.10 & OM     &  & 1.000 & 0.482 & 5.850 &   & &  &  &1.000 & 0.485 & 5.850 &  \\
 &   & & &    & RSM   &  & 0.000 & 0.000 & 1.000 &  &  &  & &  0.000 & 0.000 & 1.000  & \\
&    &  & &   & gSCAD  &   & 1.000 & 0.267 & 4.450 &  &  &  & &  1.000 & 0.000 & 3.000 &  \\
 & & &    &0.50 & OM   &  & 0.950 & 0.450 & 5.650 &   &   &  & &  1.000 & 0.465 & 5.700 &\\
 &  & & &   & RSM   &   & 0.000 & 0.000 &1.000 &  &  &   & & 0.000 & 0.000 & 1.000 & \\
  &  &  & &   & gSCAD   &  & 0.950 & 0.302 & 4.550 &   &   & &  &0.900 & 0.147 & 3.640 & \\
\hline\hline
\end{tabular}
%\end{minipage}
\end{table}
\end{center}

\bigskip

\begin{center}
\begin{table}
\caption{Results of our method    in term of  the  coverage probability (CVP), false discovery rate (FDR) and model size (MSIZE) across $200$ replications for Example 3.}\label{tab3}
\medskip
\begin{tabular}{cccccccccccccccccc}%
\hline\hline
    &  & &        &        &   &  &\multicolumn{3}{c }{Fourier basis}&     &  & &     & \multicolumn{3}{c }{B-spline basis}&  \\
\cline{8-10}\cline{15-17}
  & n & & & $\sigma$ &    &   & CVP  & FDR & MSIZE  &  &   &  & &     CVP & FDR & MSIZE  & \\
\hline
& 50 & & &0.10 &       &  & 0.80 & 0.37 & 3.20 &   & &  &  & 0.78 & 0.39 & 2.86 &  \\
 & & &    &0.25 &     &  & 0.77 & 0.29 &3.10 &   &   &  & &  0.67 & 0.43 & 3.27 &\\
 & & &    &0.50 &     &  & 0.74 & 0.34 & 3.20 &   &   &  & &  0.73 & 0.34 & 3.11 &\\
\hline
& 75 & & &0.10 &      &  & 0.65 & 0.33 & 3.75 &   & &  &  & 0.68 & 0.41 & 3.60 &  \\
 & & &    &0.25 &     &  & 0.61 & 0.34 & 3.25 &   &   &  & &  0.59 & 0.44 & 4.10 &\\
 & & &    &0.50 &     &  & 0.63 & 0.36 & 3.35 &   &   &  & &  0.55 & 0.45 & 3.90 &\\
 \hline
 &100 & & &0.10 &     &   & 0.62 & 0.30 & 4.00  &   &   & &  & 0.59& 0.46 & 3.50 &\\
  & & &    &0.25 &     &  & 0.60 & 0.30 & 4.00 &   &   &  & &  0.57 & 0.45 & 2.85 &\\
 & & &    &0.50 &     &  & 0.64 & 0.30 & 3.90 &   &   &  & &  0.57 & 0.37 & 3.10 &\\
\hline\hline
\end{tabular}
%\end{minipage}
\end{table}
\end{center}

\bigskip

\bigskip
\begin{figure}
	
	\begin{minipage}[h]{1.00\linewidth}
		\centering
		\includegraphics[width=0.6\linewidth]{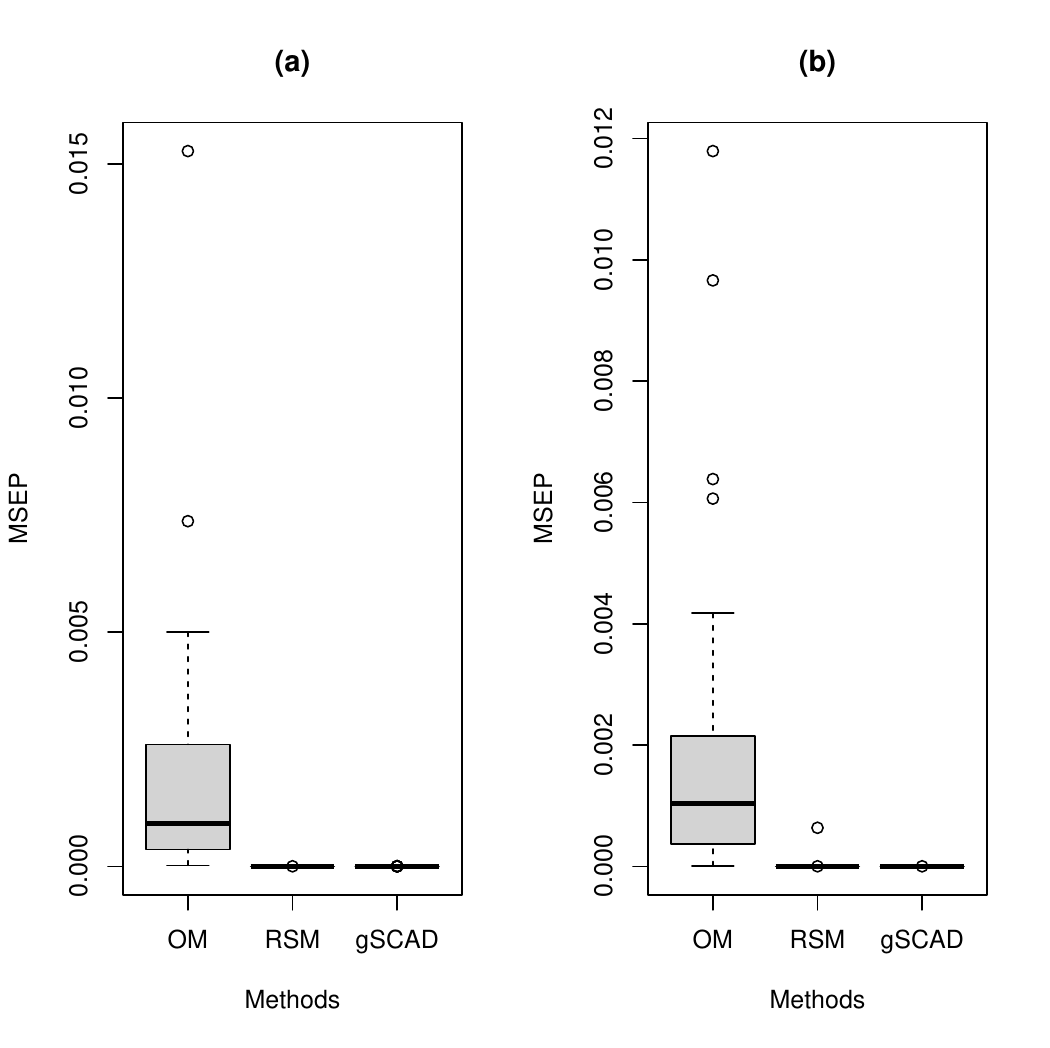}
                     \includegraphics[width=0.6\linewidth]{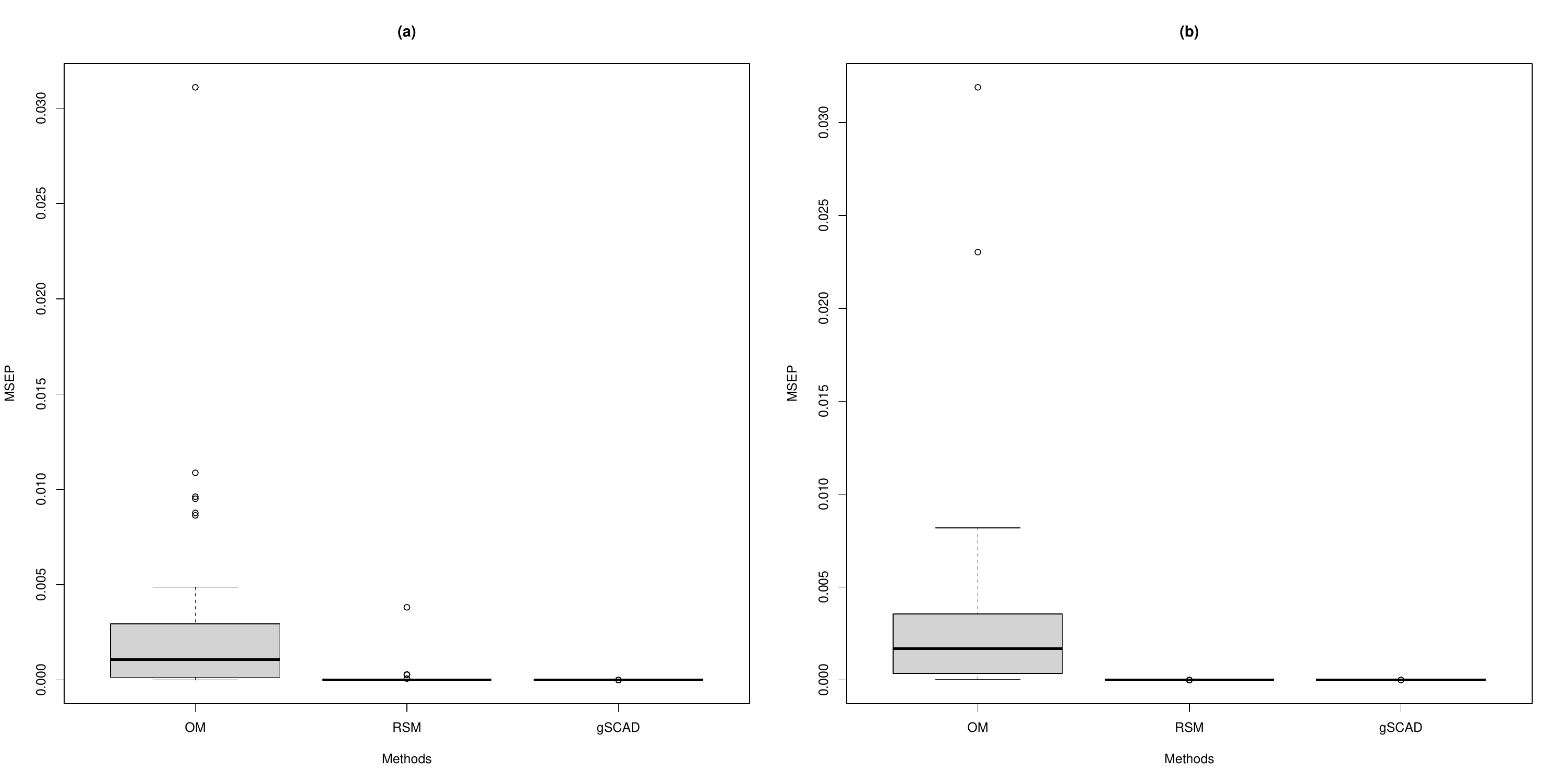}
		\caption{Boxplots showing MSEP from our method (OM), random subspace method (RSM) and group SCAD method (gSCAD) across $200$ replications for Example 1  with $n=50$, using (a) Fourier basis, (b) B-spline basis. At the top: $\sigma=0.1$; at the bottom: $\sigma=0.25$. }
	\end{minipage}
	\label{figure1}
	%\hfill
	
\end{figure}

\begin{figure}
	
	\begin{minipage}[h]{1.00\linewidth}
		\centering
		\includegraphics[width=0.6\linewidth]{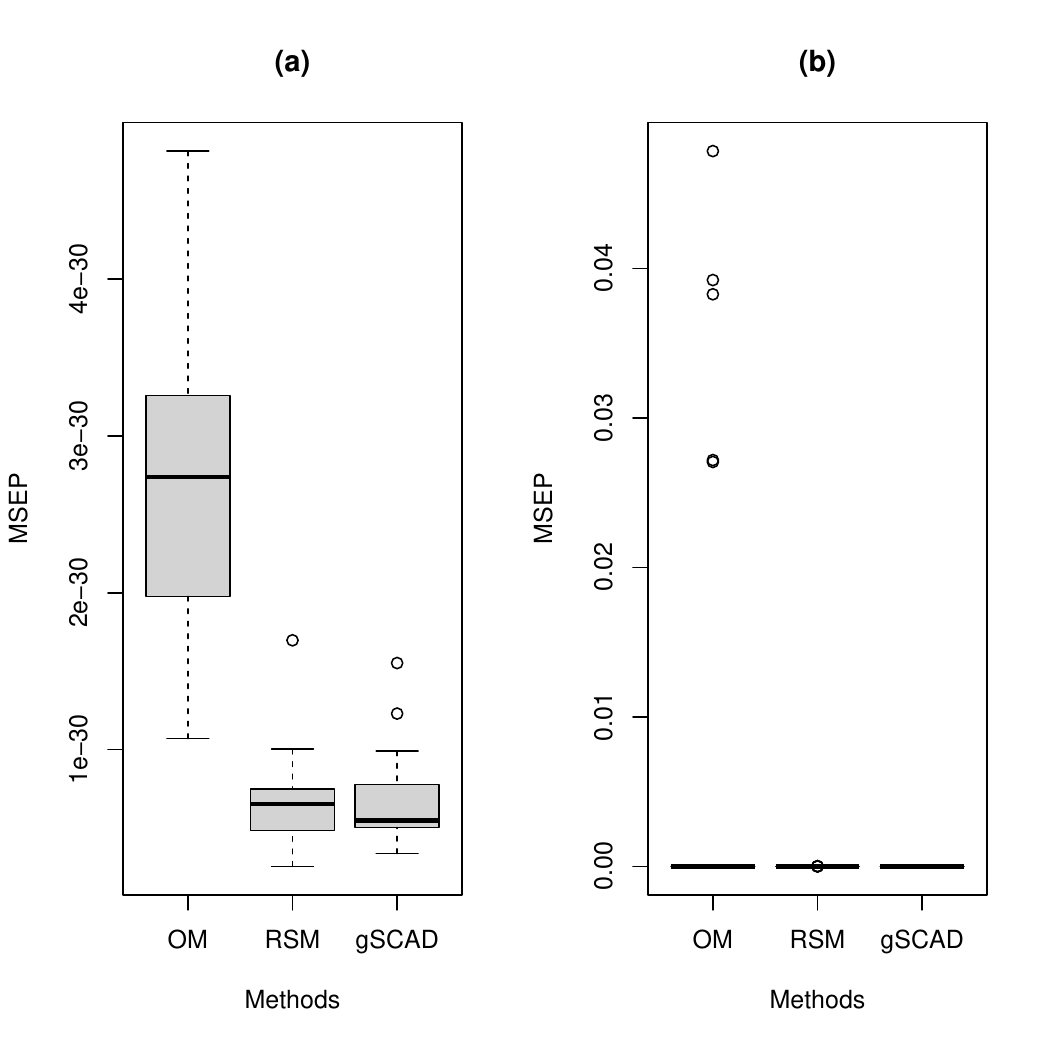}
                      \includegraphics[width=0.6\linewidth]{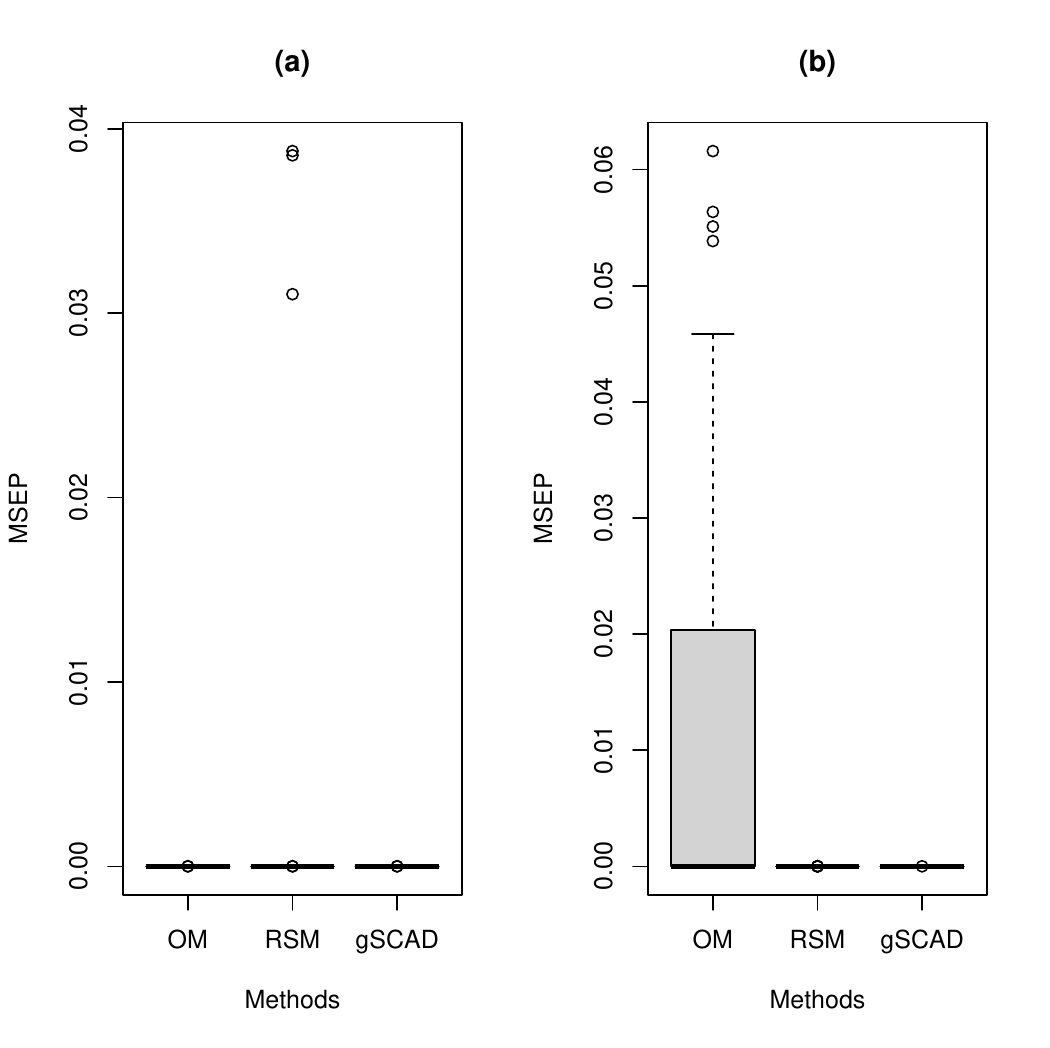}
		\caption{Boxplots showing MSEP from our method (OM), random subspace method (RSM) and group SCAD method (gSCAD) across $200$ replications for Example 1  with $n=75$, using (a) Fourier basis, (b) B-spline basis. At the top: $\sigma=0.1$; at the bottom: $\sigma=0.25$.}
	\end{minipage}
	\label{figure2}
	%\hfill
	
\end{figure}

\begin{figure}
	
	\begin{minipage}[h]{1.00\linewidth}
		\centering
		\includegraphics[width=0.55\linewidth]{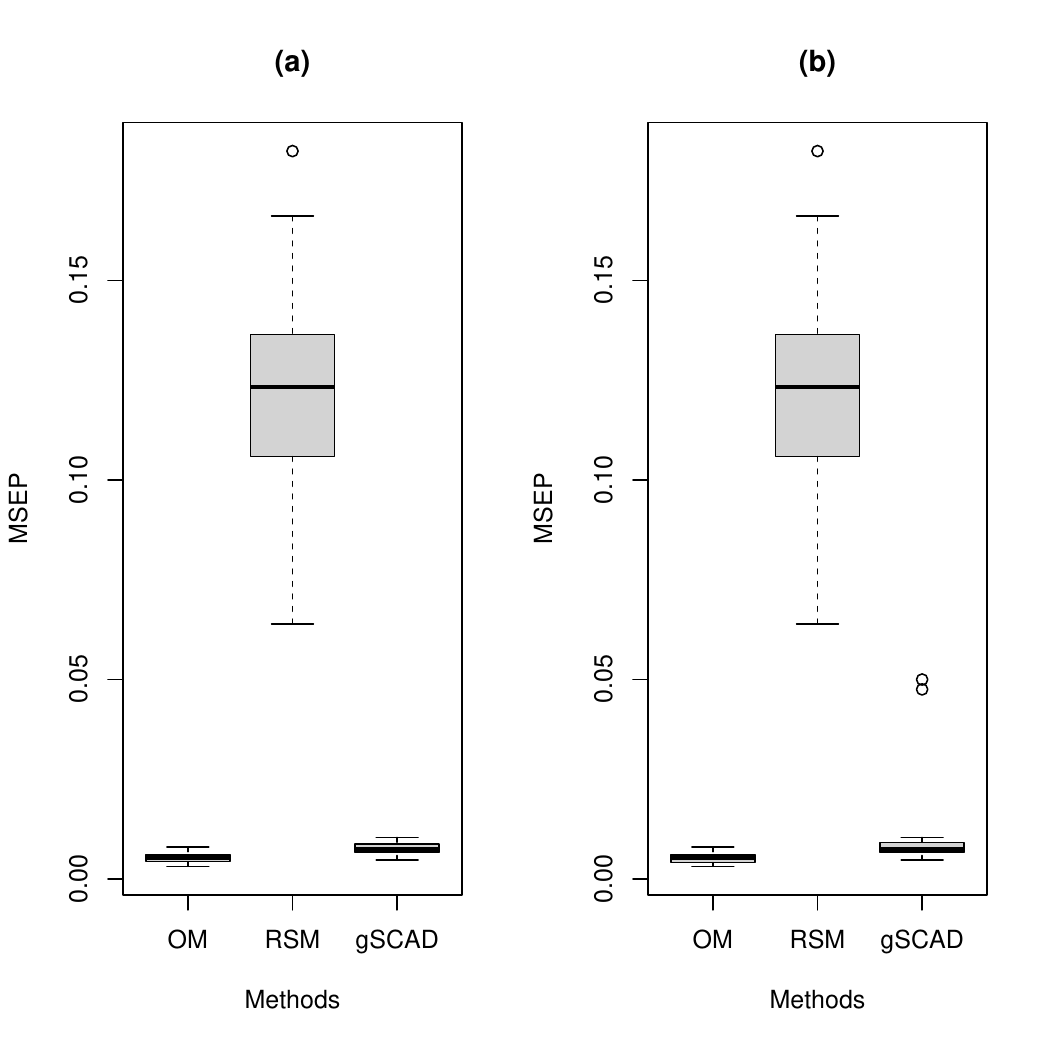}
                    \includegraphics[width=0.55\linewidth]{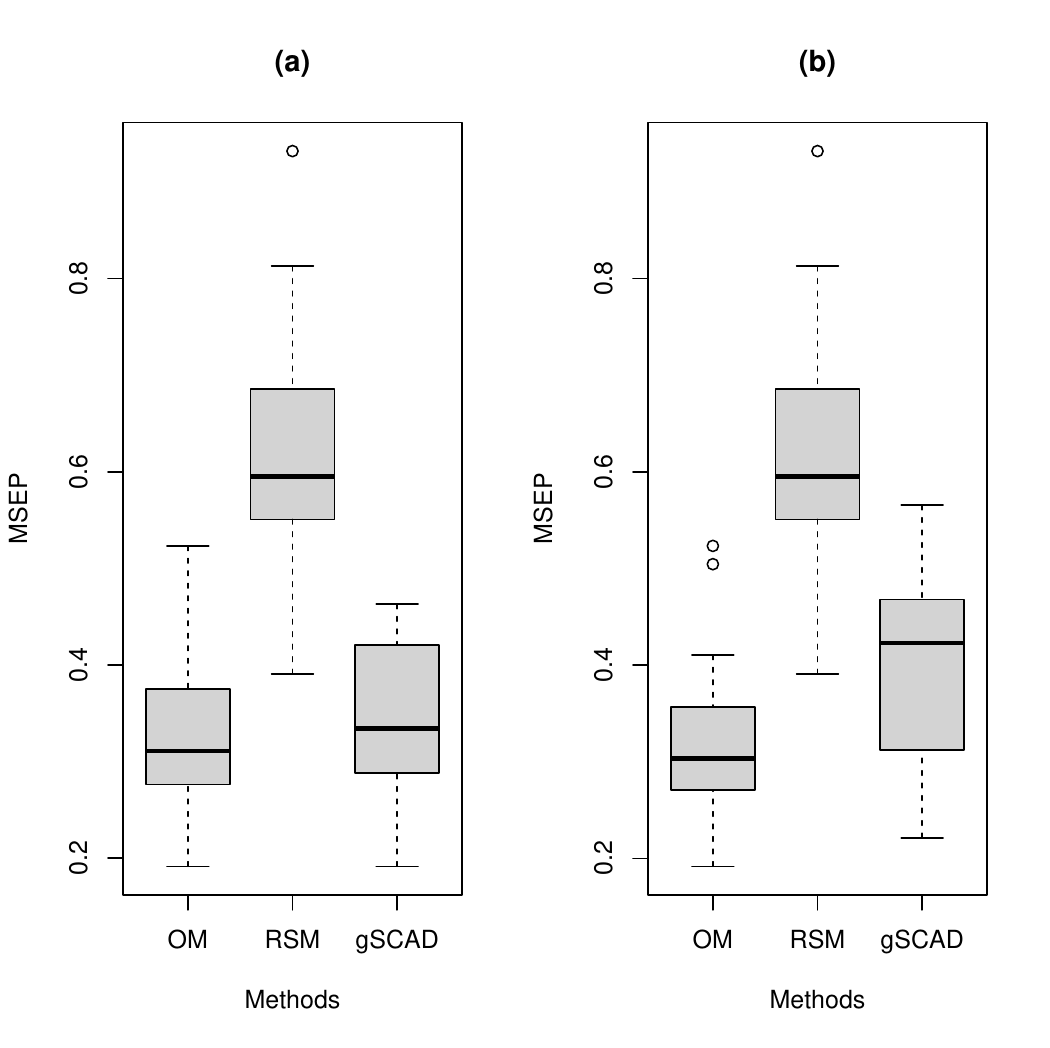}
		\caption{Boxplots showing MSEP from our method (OM), random subspace method (RSM) and group SCAD method (gSCAD) across $200$ replications for Example 2  with $n=50$, using (a) Fourier basis, (b) B-spline basis. At the top: $\sigma=0.1$;    at the bottom: $\sigma=0.5$.}
	\end{minipage}
	\label{figure3}
	%\hfill
	
\end{figure}

\begin{figure}
	
	\begin{minipage}[h]{1.00\linewidth}
		\centering
		\includegraphics[width=0.6\linewidth]{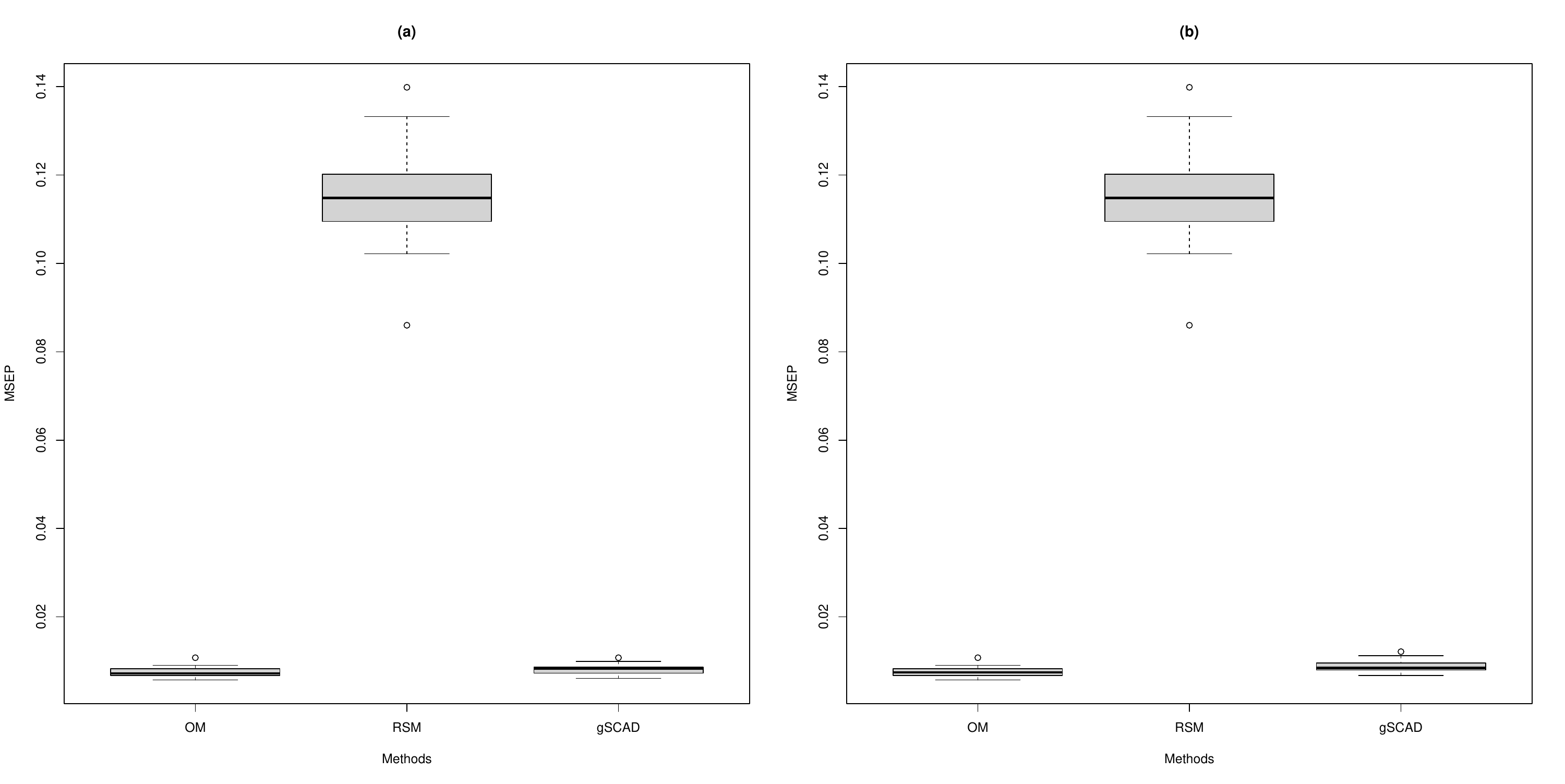}
                    \includegraphics[width=0.6\linewidth]{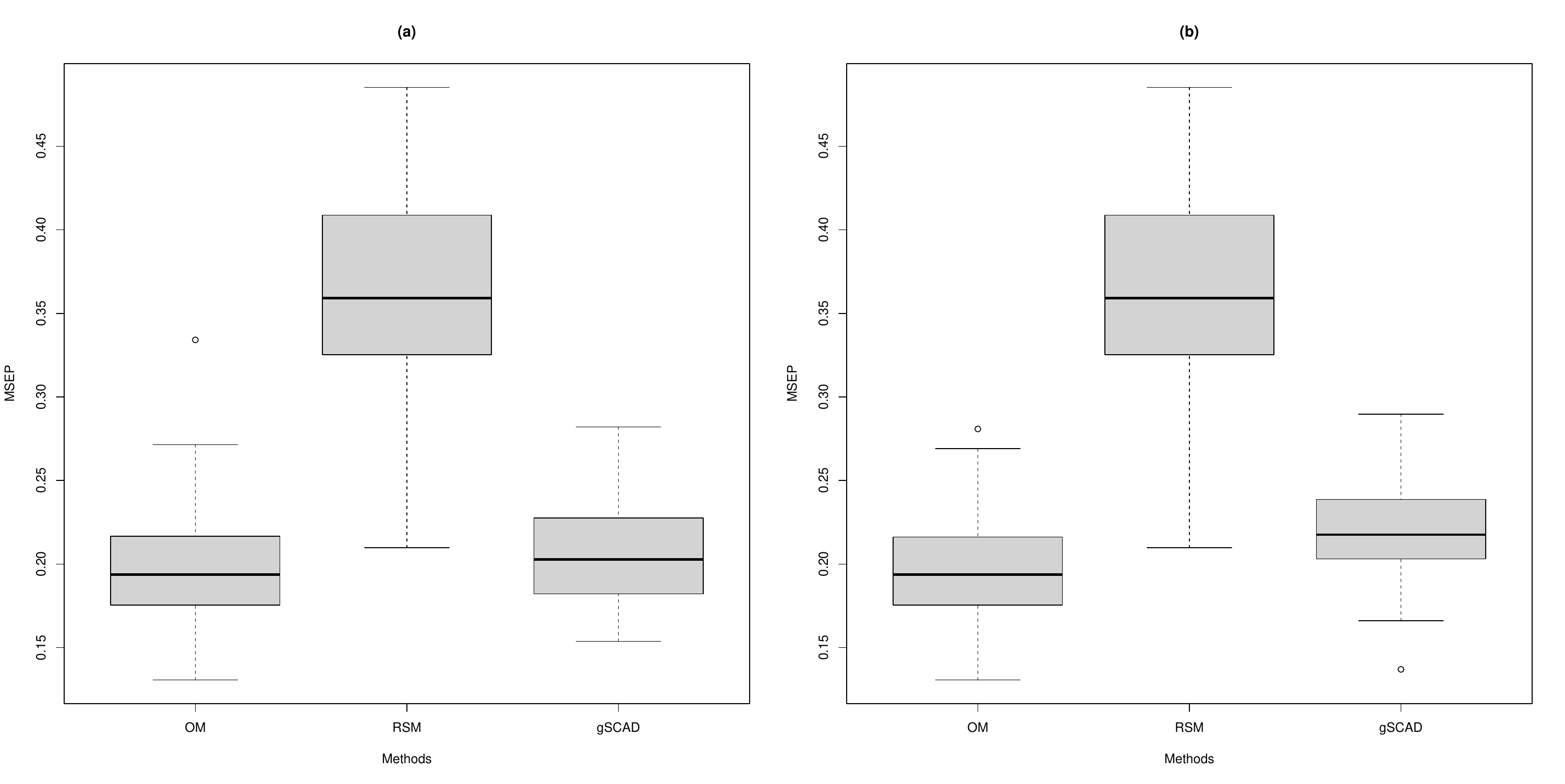}
		\caption{Boxplots showing MSEP from our method (OM), random subspace method (RSM) and group SCAD method (gSCAD) across $200$ replications for Example 2  with $n=100$, using (a) Fourier basis, (b) B-spline basis. At the top: $\sigma=0.1$;  at the bottom: $\sigma=0.5$.}
	\end{minipage}
	\label{figure4}
	%\hfill
	
\end{figure}

%**********

\begin{figure}
	
	\begin{minipage}[h]{1.00\linewidth}
		\centering
		\includegraphics[width=0.5\linewidth]{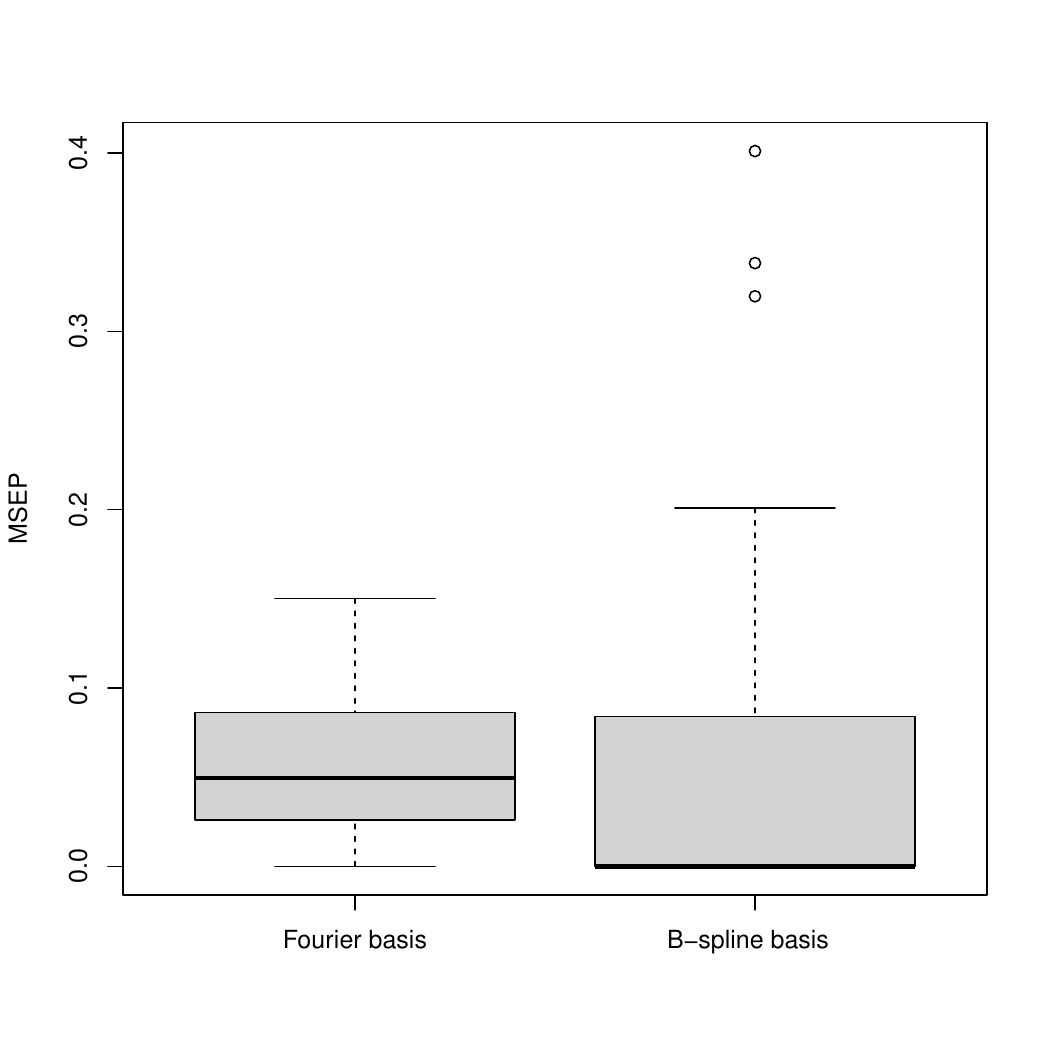}
                    \includegraphics[width=0.5\linewidth]{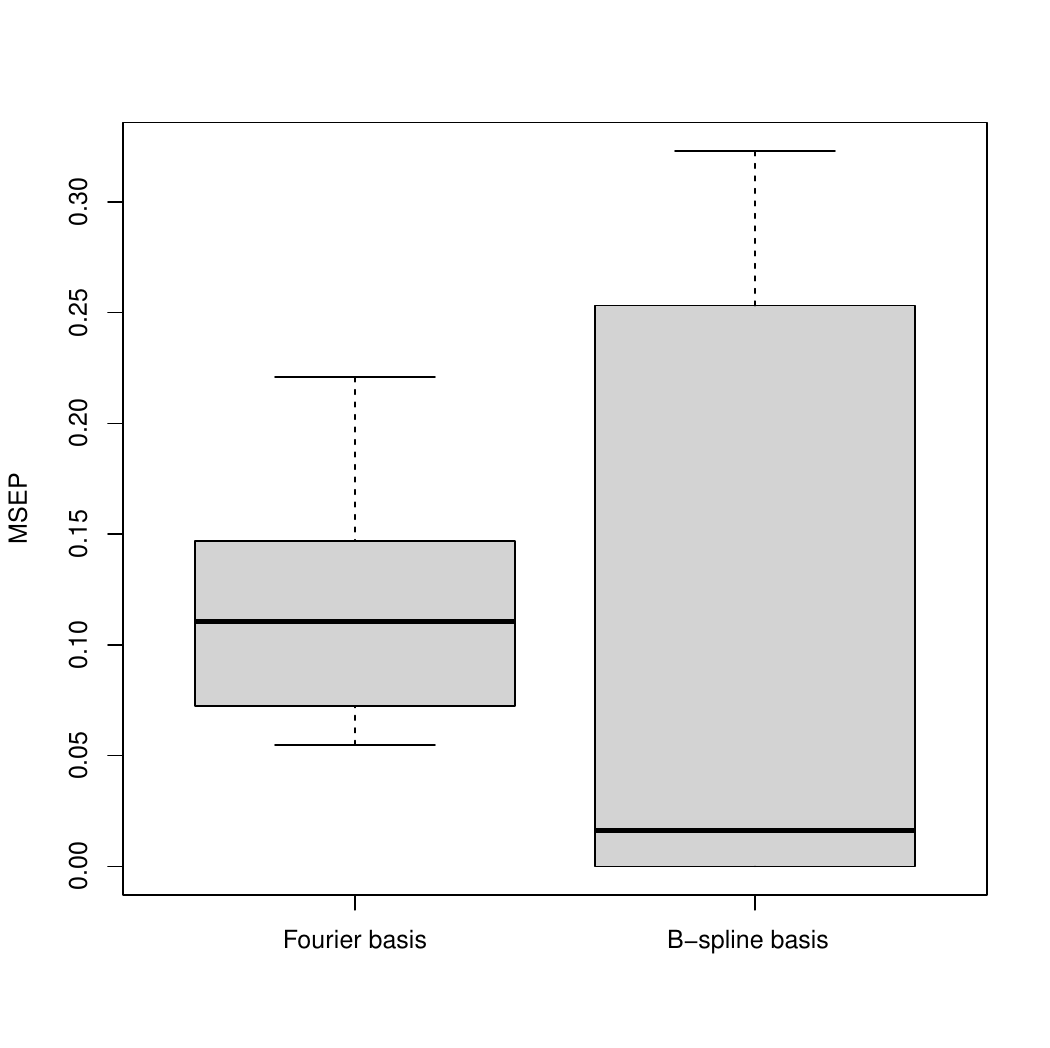}
		\caption{Boxplots showing MSEP from our method  across $200$ replications for Example 3  with $n=50$. At the top: $\sigma=0.1$;   at the bottom: $\sigma=0.5$.}
	\end{minipage}
	\label{figure5}
	%\hfill
	
\end{figure}

\begin{figure}
	
	\begin{minipage}[h]{1.00\linewidth}
		\centering
		\includegraphics[width=0.5\linewidth]{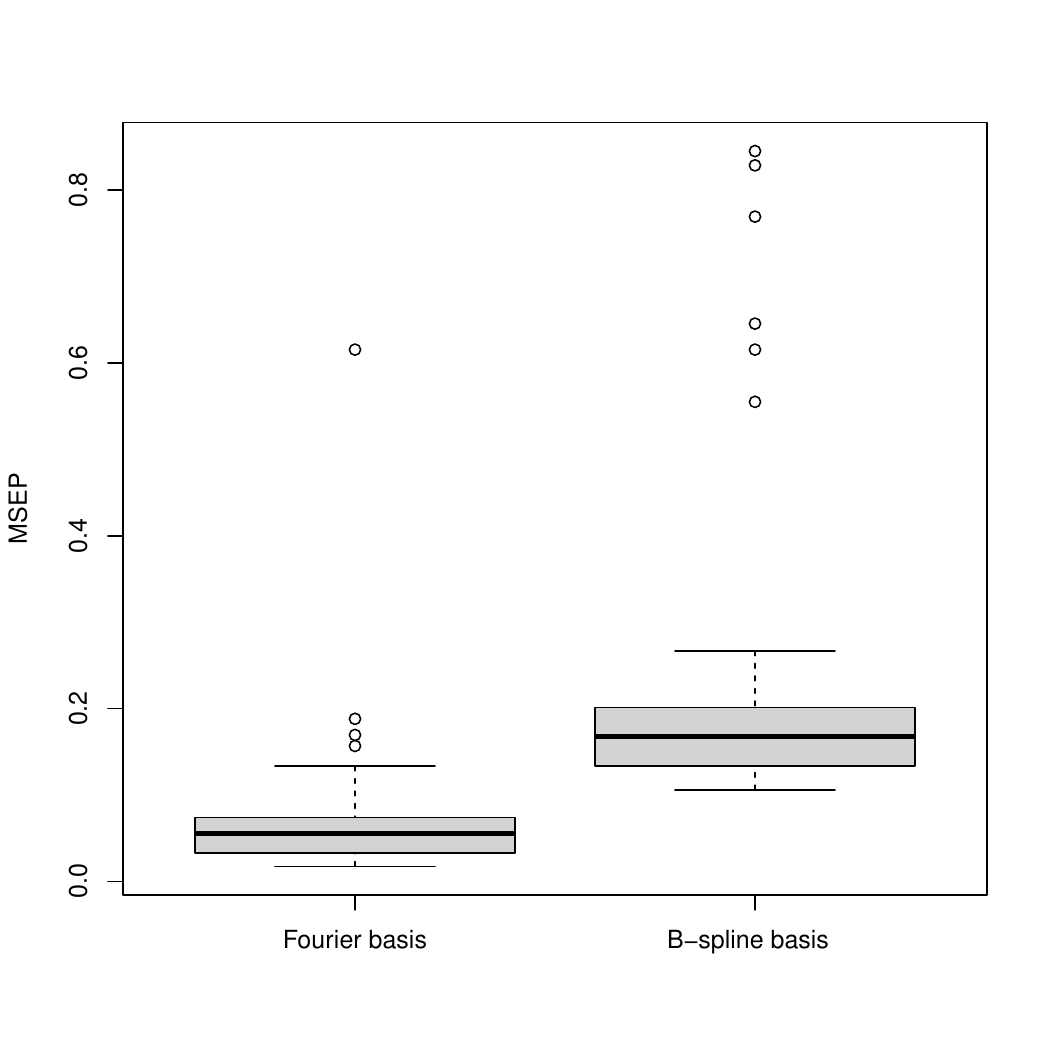}
                    \includegraphics[width=0.5\linewidth]{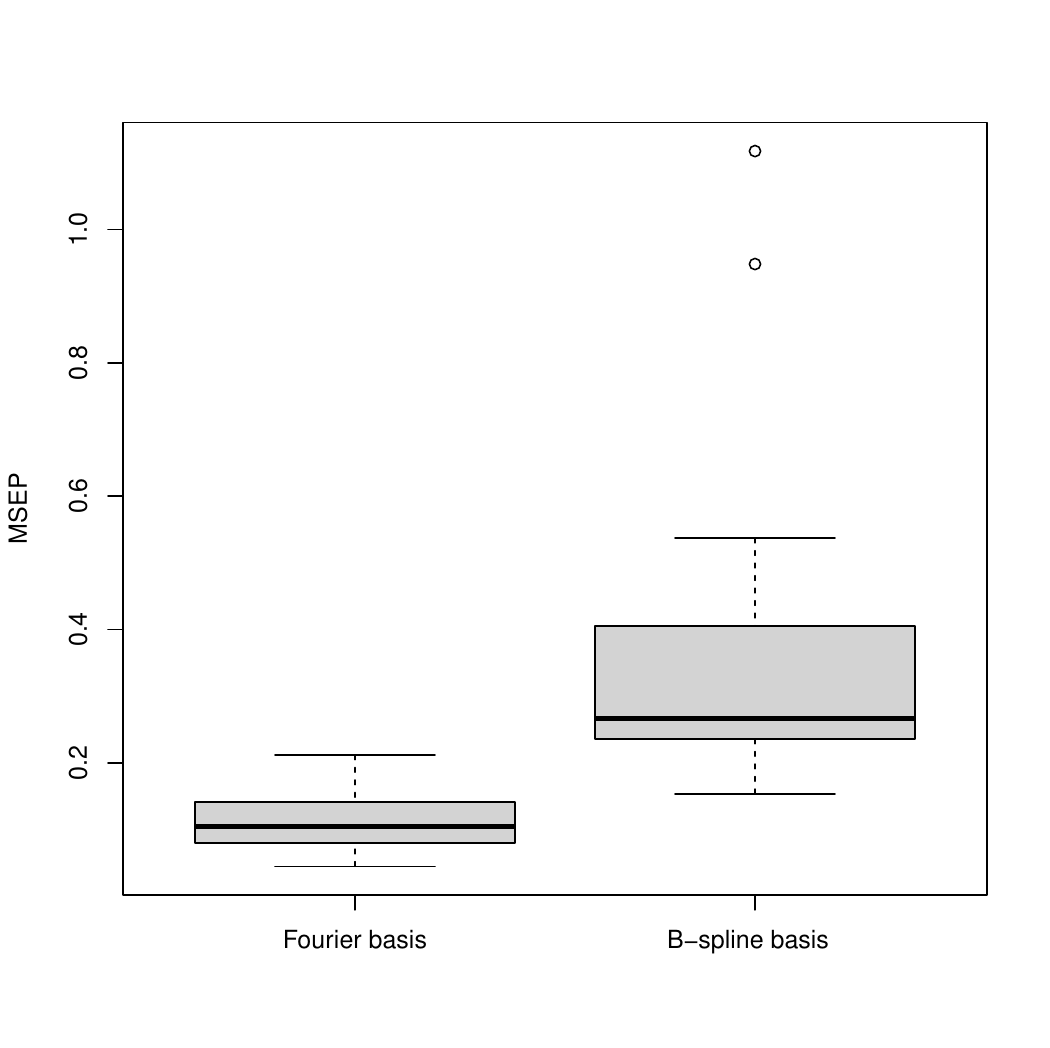}
		\caption{Boxplots showing MSEP from our method  across $200$ replications for Example 3  with $n=100$. At the top: $\sigma=0.1$; at the bottom: $\sigma=0.5$.}
	\end{minipage}
	\label{figure6}
	%\hfill
	
\end{figure}

\section{Conclusion}
\label{sec:conc}
\noindent We have proposed  a new approach for variable selection in a multivariate functional linear regression model, by extending to this framework a method that has already been used in multivariate linear regression.  This extension was  made possible thanks to  basis representations of the functional explanatory variables. One of the advantages of this method is that it can be used in the functional multivariate case, which is not the case for most existing methods for  variable selection  in functional regression models  which only deal with the univariate case.  Numerical experiments  have shown its good performance in comparison with  the random subspace method   of Smaga and Matsui (2018)  and the group SCAD method of  Matsui and Konishi (2001). Thus, we have provided a new competitive
alternative to perform functional variable selection.

\bigskip

\medskip

%\section*{References}

%\section*{References}

% To ensure accuracy, get them from MathSciNet whenever possible. Typeset them with BibTeX using JMVA's style file, \texttt{myjmva.bst}.
%\bibliography{}
%\bibliographystyle{myjmva}
%\begin{thebibliography}
%\bibliography{trial}
%\end{document}

%\bibliographystyle{myjmva}
%\section*{}

\end{document}